\documentclass[preprint]{elsarticle}
\usepackage{graphicx}
\usepackage{amssymb,amsmath}
\usepackage{subfig}
\newcommand{\nc}{\newcommand}
\nc{\rnc}{\renewcommand}
\nc{\x}[1]{\mbox{#1}}           
\nc{\hs}[1]{\hspace*{#1}}
\nc{\vs}[1]{\vspace*{#1}}
\rnc{\theta}{\vartheta}
\rnc{\rho}{\varrho}
\rnc{\epsilon}{\varepsilon}
\nc{\dd}{{\mathrm{d}}}
\nc{\ii}{{\mathrm{i}}}
\nc{\ee}{{\mathrm{e}}}
\nc{\mm}[1]{\mathbf{#1}}
\nc{\mms}[1]{\boldsymbol{#1}}
\nc{\suml}[2]{\sum \limits_{#1}^{#2}}
\nc{\intl}[2]{\int \limits_{#1}^{#2}}
\rnc{\matrix}[2]{\left[\!\!\begin{array}{#1} #2\end{array}\!\!\right]}
\rnc{\vector}[1]{\matrix{c}{#1}}
\nc{\inv}{^{-1}}
\nc{\tra}{^{\mathrm T}}
\nc{\sgn}{\mathrm{sgn}}
\nc{\diag}{\mbox{\bf diag}}
\nc{\eo}{EO}
\nc{\nh}{{N_{\mathrm h}}}
\nc{\ncp}{{N_{\mathrm{CP}}}}
\nc{\nnl}{{N_{\mathrm{N}}}}
\nc{\nll}{{N_{\mathrm{L}}}}
\nc{\ndim}{{N_{\mathrm{dof}}}}
\nc{\nnlred}{{N_{\mathrm{nl,red}}}}
\nc{\myquote}[1]{`#1'}
\nc{\kreis}[1]{\mbox{\mbox{\text{\Large{\ensuremath{\bigcirc}}}}\hspace{-2.4ex}#1\hspace{1.2ex}}}
\nc{\prob}{\operatorname{Prob}}
\nc{\fnp}{N}
\nc{\fnpopt}{N^{\rm{opt}}}
\nc{\fnpoptc}{N^{\rm{opt}}_{\rm{c}}}
\nc{\fnpoptuc}{N^{\rm{opt}}_{\rm{uc}}}
\nc{\ares}{a^{\rm{res}}}
\nc{\astick}{a^{\rm{res}}_{\rm{stick}}}
\nc{\amin}{a^{\mathrm{opt}}_{\mathrm{c}}}
\nc{\fres}{f^{\rm{res}}}
\nc{\areference}{a^{\rm{ref}}}
\nc{\argm}{\operatorname{arg}}
\nc{\kn}{{k_{\mathrm n}}}
\nc{\kt}{{k_{\mathrm t}}}
\rnc{\Re}[1]{\operatorname{Re}\lbrace #1 \rbrace}
\nc{\prog}[1]{{\sf{#1}}}
\nc{\matlab}{\prog{Matlab}}
\nc{\ie}{i.\,e.\,}
\nc{\eg}{e.\,g.\,}
\nc{\cf}{cf.\,}
\nc{\etal}{et~al.\,}
\nc{\z}[2]{\x{#1}~{\x{\cite{#2}}}}
\nc{\zo}[1]{\x{\cite{#1}}}
\nc{\fabstand}{\,}
\nc{\fp}{\fabstand.}
\nc{\fk}{\fabstand,}
\nc{\tab}[5][tbh]{\begin{table}[#1]\centering\caption{#4\label{tab:#5}}\begin{tabular}{#2}\hline #3 \\ \hline\end{tabular}\end{table}}
\newcommand{\fss}[4][tbh]{%
 \begin{figure}[#1]
 \centering
 \if@draft
  \framebox[150mm]{\raisebox{0mm}[25mm][25mm]{\texttt{#2}}}
 \else
  \includegraphics[scale=#4]{#2}
 \fi
 \caption{#3}
 \label{fig:#2}
\end{figure}
}
\newcommand{\myf}[8][tbh]{%
    \begin{figure}[#1]
        \centering
        \subfloat[#4 \label{fig:#2}]{
            \includegraphics[width=#6\textwidth]{#2}
        }
        \hspace{0.25cm}
        \subfloat[#5 \label{fig:#3}]{
            \includegraphics[width=#7\textwidth]{#3}
        }
        \caption{#8}
    \end{figure}
}
\newcommand{\mythreef}[8][tbh]{%
    \begin{figure}[#1]
        \centering
        \subfloat[#5 \label{fig:#2}]{
            \includegraphics[width=0.29\textwidth]{#2}
        }
        \hspace{0.15cm}
        \subfloat[#6 \label{fig:#3}]{
            \includegraphics[width=0.29\textwidth]{#3}
        }
        \hspace{0.15cm}
        \subfloat[#7 \label{fig:#4}]{
            \includegraphics[width=0.29\textwidth]{#4}
        }
        \caption{#8}
    \end{figure}
}
\nc{\e}[2]{\begin{equation} #1 \label {eq:#2} \end{equation}}
\nc{\ea}[2]{
\begin{eqnarray}
#1 \label {eq:#2} \end{eqnarray}}
\nc{\eal}[3][0.0ex]{
\begin{samepage}
\begin{eqnarray*}
#2
\end{eqnarray*}
\nopagebreak[4] \vs{#1} \nopagebreak[4]\vs{-2ex} \nopagebreak[4]
\begin{eqnarray}
\label {eq:#3}
\end{eqnarray}
\end{samepage}\hs{-0.35em}}
\nc{\g}[1]{{$#1$}}
\nc{\gehg}[3]{\x{${#1}= #2 ~ {\rm{#3}}$}}
\nc{\fref}[1]{{Fig.~\ref{fig:#1}}}
\nc{\frefo}[1]{{\ref{fig:#1}}}
\nc{\frefoo}[1]{{#1}}
\nc{\frefs}[1]{{Figs.~\ref{fig:#1}}}
\nc{\tref}[1]{{Tab.~\ref{tab:#1}}}
\nc{\trefo}[1]{{\ref{tab:#1}}}
\nc{\trefs}[1]{{Tab.~\ref{tab:#1}}}
\nc{\erefn}[1]{{Eq.~(\ref{#1})}}
\nc{\eref}[1]{{Eq.~(\ref{eq:#1})}}
\nc{\erefo}[1]{(\ref{eq:#1})}
\nc{\erefs}[2]{{Eqs.~(\ref{eq:#1}) and (\ref{eq:#2})}}
\nc{\erefss}[1]{{Eqs.~(\ref{eq:#1})}}
\nc{\sref}[1]{{Section~\ref{sec:#1}}}
\nc{\srefo}[1]{\ref{sec:#1}}
\nc{\srefs}[1]{{Sections~\ref{sec:#1}}}
\nc{\aref}[1]{{{Appendix~\ref{asec:#1}}}}
\nc{\arefo}[1]{{\ref{asec:#1}}}
\nc{\arefs}[1]{{{Appendices~\ref{asec:#1}}}}
\nc{\ssref}[1]{{Subsection~\ref{sec:#1}}}
\nc{\ssrefo}[1]{\ref{sec:#1}}
\nc{\ssrefs}[1]{{Subsections~\ref{sec:#1}}}
\begin{document}

\begin{frontmatter}

\title{A METHOD FOR NONLINEAR MODAL ANALYSIS AND SYNTHESIS: APPLICATION TO HARMONICALLY FORCED AND
SELF-EXCITED MECHANICAL SYSTEMS}

\author[ids]{Malte Krack\corref{cor1}}
\ead{krack@ila.uni-stuttgart.de}
\author[ids]{Lars Panning-von Scheidt}
\author[ids]{J\"org Wallaschek}


\cortext[cor1]{Corresponding author}

\begin{abstract}
\textit{The recently developed generalized Fourier-Galerkin method
is complemented by a numerical continuation with respect to the
kinetic energy, which extends the framework to the investigation of
modal interactions resulting in folds of the nonlinear modes. In
order to enhance the practicability regarding the investigation of
complex large-scale systems, it is proposed to provide analytical
gradients and exploit sparsity of the nonlinear part of the
governing algebraic equations.\\
A novel reduced order model (ROM) is developed for those regimes
where internal resonances are absent. The approach allows for an
accurate approximation of the multi-harmonic content of the resonant
mode and accounts for the contributions of the off-resonant modes in
their linearized forms. The ROM facilitates the efficient analysis
of self-excited limit cycle oscillations, frequency response
functions and the direct tracing of forced resonances. The ROM is
equipped with a large parameter space including parameters
associated with linear damping and near-resonant harmonic forcing
terms. An important objective of this paper is to demonstrate the
broad applicability of the proposed overall methodology. This is
achieved by selected numerical examples including finite element
models of structures with strongly nonlinear, non-conservative
contact constraints.}
\end{abstract}


\end{frontmatter}

\section{Introduction\label{sec:introduction}}
Akin to its linear counterpart, nonlinear modal analysis is
particularly suited for the analysis of dynamical systems. Modal
analysis facilitates understanding of the energy-dependent system
behavior in nonlinear systems regarding eigenfrequencies, modal
damping, stiffening/softening characteristics, localization effects
and internal resonances. The concept of nonlinear modes dates back
to \z{Rosenberg}{rose1960} and the interested reader is referred to
\zo{vaka2008,kers2009} for a good overview on various concepts and
theories. Despite the fact that superposition and orthogonality
conditions are not valid in the nonlinear case, nonlinear modes have
been widely used for the approximate synthesis of forced vibrations
\zo{szem1979,jeze1991b,chon2000a,gibe2003,laxa2009}.\\
In spite of their opportunities for qualitative and quantitative
analysis of nonlinear systems, methods related to nonlinear modes
are seldom applied to industrial problems. In the authors' opinion,
the reasons are that most methods are restricted to smooth and
conservative nonlinearities and rarely proved to cope with systems
featuring many degrees of freedom (DOFs) such as large scale finite
element models of typical industrial applications.\\
Several methods have been developed for analytical and numerical
calculation of nonlinear modes in the past. We will focus on those
methods which are well-suited for systems with generic and strong
nonlinearities. Perturbation techniques such as the normal form
approach \zo{jeze1991} and the method of multiple scales
\zo{nayf1979} are not considered since they are restricted to small
degree polynomial nonlinearities.\\
A method of broad applicability is the invariant manifold approach,
as proposed in \zo{shaw1993,nayf2000,jian2005,pierre2006,touz2006}.
It is based on the invariance property of certain periodic motions
of the system, \ie a nonlinear mode is defined as an invariant
relationship (manifold) between several master coordinates and the
remaining coordinates of the system. This manifold can be governed
by partial differential equations arising from the substitution of
the manifold into the state form of the equations of motion. For the
solution of the governing equations, asymptotic expansions were
originally employed \zo{shaw1993,touz2006} and later a more general
Galerkin ansatz was developed \zo{pesh2002} to increase the accuracy
of this approach. The invariant manifold approach was extended to
account for the effect of harmonic excitation \zo{jian2005} and
viscous damping \zo{touz2006}. Various systems have been studied,
including piecewise linear systems \zo{jian2004} and systems with
internally resonant nonlinear modes \zo{pierre2006}. However, in
these investigations, the focus of application was clearly set on
small-scale systems with conservative
nonlinearities.\\
Another class of methods for the determination of the modal
properties of nonlinear systems is based on parameter
identification, see \eg \zo{chon2000a,gibe2003,kers2005}. Response
data, obtained either by simulation or measurement, is gathered and
modal parameters are identified by fitting original response data to
data from nonlinear modal synthesis. An important subclass of
parameter identification methods is the so called force
appropriation, where the objective is to extract modal properties
from the resonance reached by suitably adjusting the forcing
parameters \zo{zasp2007,peet2010,kuet2012}. The weak point of these
strategies is clearly their signal-dependent nature and the fact
that the modal parameters are typically extracted from a forced
rather than from an autonomous system. Moreover, further numerical
or experimental effort is required to obtain the response data. One
of the main benefits of this method is that typically no model is
required for the nonlinearities which enables broad applicability.\\
More recently, \z{Kerschen \etal}{kers2009,peet2009} developed a
method for the calculation of nonlinear modes of conservative
large-scale mechanical systems. Their method exploits the periodic
nature of nonlinear normal modes and is based on the shooting
algorithm in conjunction with time-step integration. Unfortunately,
it is not yet clear, whether this method can be extended to
dissipative systems.\\
The periodicity of nonlinear modes in conservative systems is also
the starting point for the application of the Harmonic Balance
Method (HBM) \zo{leun1992,ribe2000,coch2009}. The HBM is well-suited
for the analysis of strongly nonlinear systems and often leads to
reduced computational effort compared to time integration
approaches. \z{Laxalde \etal}{laxa2008a,laxa2009} generalized the
HBM to account for the energy decay of the nonlinear mode in
dissipative systems. Their modal analysis technique is therefore
qualified for the investigation of non-conservative nonlinear
systems. Modal properties have also been exploited for the forced
response synthesis and the computation of limit-cycle-oscillations.
The authors applied the methodology to turbomachinery bladed disks
with friction interfaces featuring constant normal load.\\
The goal of the present paper is twofold. Firstly, a methodology is
developed for efficient numerical computation of nonlinear modes of
large-scale mechanical systems with generic, including strong and
non-smooth, conservative or dissipative nonlinearities, see
\sref{cnma}. Secondly, nonlinear modal properties are used to
accurately calculate forced and self-excited vibrations in
\sref{nlrom}. In \sref{numerical_examples}, the proposed methodology
is applied to several nonlinear mechanical systems including systems
with friction and unilateral contact, and strengths and weaknesses
compared to conventional methods are discussed.


\section{Complex Nonlinear Modal Analysis\label{sec:cnma}}
The equations of motion of a discrete, time-invariant, autonomous
mechanical system can be stated as
\e{\mm{M}\ddot{\mm{u}}(t)+ \mm{K}\mm{u}(t)+
\mm{g}\left(\mm{u}(t),\dot{\mm{u}}(t)\right)
=\mm{0}\fp}{eqm_autonomous}
Herein, \g{\mm M = \mm M\tra>0} is the real, symmetric, positive
definite mass matrix, \g{\mm K=\mm K\tra} is the real, symmetric
stiffness matrix and \g{\mm u(t)} is the vector of generalized
coordinates. The vector \g{\mm g} can comprise linear and nonlinear
functions dependent on displacement and velocity. Without loss of
generality, the generalized coordinates of the system can be defined
in such a way that \g{\mm u(t)=\mm 0} is an equilibrium point and
\g{\mm K} contains the symmetric part of the linearization of \g{\mm
g} with respect to \g{\mm u} around this equilibrium. The number of
DOFs is denoted \g{\ndim}. The restrictions made regarding symmetry
of the structure in \eref{eqm_autonomous} are relaxed in
\ssref{ana_general}. It should be emphasized that the use of
generalized coordinates in \eref{eqm_autonomous} explicitly allows
preceding component mode synthesis, which can be very useful when
treating large-scale structures with localized
nonlinearities.\\
Non-trivial solutions \g{\mm u(t)} of \eref{eqm_autonomous} are
sought in the form of a generalized Fourier ansatz \zo{laxa2009},
\e{\mm u(t) = \Re{\suml{n=0}{\nh}\mm U_n\ee^{n\lambda
t}}\fp}{fourier_ansatz}
Herein, \g{\lambda = -D\omega_0+\ii\omega_0\sqrt{1-D^2}} is the
complex eigenvalue with the eigenfrequency \g{\omega_0} and the
modal damping ratio \g{D}, and \g{\mm U_n} are vectors of complex
amplitudes. The ansatz in \eref{fourier_ansatz} induces the
assumption that the damping of the system is frequency-independent
\zo{laxa2009}. For conservative systems the damping is zero,
\g{D=0}, so that of conservative systems are not
affected by this assumption.\\
If only the first harmonic \g{n=1} is retained, the ansatz
degenerates to the well-known exponential ansatz for damped linear
systems. Further, for \g{D=0} \eref{fourier_ansatz} is completely
equivalent to the conventional HBM ansatz for conservative
autonomous systems. The damping term \g{D} takes into account the
energy decay
of the nonlinear mode.\\
Inserting ansatz \erefo{fourier_ansatz} into \eref{eqm_autonomous}
and subsequent Fourier-Galerkin projection with respect to the base
functions gives rise to a system of nonlinear algebraic equations,
\e{\mm S_n(\lambda)\mm U_n + \mm G_n\left(\mm U_0\,,\,\cdots\,,\mm
U_{\nh}\right) = \mm 0\,,\quad n = 0,\cdots,\nh\fp}{complex_evp}
Capital letters \g{\mm U_n\,,\mm G_n} in this equation denote
complex amplitudes of the corresponding lower-case time-domain
variables. \g{\mm S_n} are the blocks of the dynamic stiffness
matrix,
\e{\mm S_n = (n\lambda)^2\mm M + \mm K\fp}{dyn_stiff}

\subsection{Mode normalization\label{sec:ana_norm}}
The number of unknowns in \eref{complex_evp} exceeds the number of
equations by two. Phase and amplitude normalization have to be
performed. In \zo{laxa2009}, normalization by prescribed first
complex amplitude \g{q_{\rm m}} and phase \g{\phi_{\rm m}} of a
specified coordinate \myquote{m} was proposed,
\e{\left| U_1^{(\mathrm m)} \right| - q_{\mathrm m} = 0\,\,\land
\,\, \argm\left(U_1^{(\mathrm m)}\right) - \phi_{\rm m} =
0\fp}{amplitude_normalization}
An amplitude normalization with respect to the kinetic energy
facilitates a direct calculation of the
frequency-energy-relationship. This is required to resolve modal
interactions \zo{lee2005,kers2009}. The corresponding normalization
conditions thus read
\e{\frac{1}{T}\intl{0}{T}\frac12\mm{\dot u}^{\mathrm T}\mm M\mm{\dot
u}\dd t - q_{\mathrm{kin}} = 0 \,\,\land \,\,
\argm\left(U_1^{(\mathrm m)}\right) - \phi_{\rm m} =
0\fp}{energy_normalization}
For clarity, the time dependence of variables is not denoted here
and in the following. The kinetic energy is represented by its mean
value on the pseudo-period \g{T=\frac{2\pi}{\omega_0}}. Note that in
conjunction with the Fourier ansatz, this integral can be easily
evaluated using Parseval's theorem.\\
Strong local nonlinearities can induce abrupt changes of the mode
shape in the vicinity of the source of nonlinearity, \eg close to a
contact area. If a DOF in such a region is chosen for the amplitude
normalization, a weak numerical performance is possible. In
contrast, global features such as the kinetic energy typically
exhibit a smoother relationship with the modal properties.
Therefore, it is expected that the energy normalization can
generally improve the computational robustness of the analysis.\\
The choice of the phase \g{\phi_{\rm m}} is arbitrary in an
autonomous system. A practical reformulation of the generalized
phase condition is to simply set the real or the imaginary part of
the component to zero.\\
An appropriate master coordinate must be specified for the phase and
amplitude normalization. Note that if \g{U_1^{(\mathrm m)}=0},
\erefs{amplitude_normalization}{energy_normalization} do not allow a
unique normalization of the mode. However, for some types of
nonlinearities such as contact constraints, it is generally possible
that certain DOF are fully stuck at specific energy levels. In order
to avoid this particular situation, it not recommended to use a
nonlinear DOF as master coordinate for the mode normalization in
presence of contact. In \zo{nayf2000} it is proposed to specify the
amplitude of the linearized mode to be analyzed as master
coordinate.

\subsection{Evaluation of the nonlinear terms\label{sec:ana_nlforces}}
As in \zo{laxa2009}, the nonlinear terms \g{\mm G_n} are integrated
on the pseudo-period and therefore the coordinates and forces are
treated as periodic within this step, in contrast to the ansatz
given by \eref{fourier_ansatz}, which takes into account the energy
decay. One advantage of this strategy is that classical HBM
frameworks do not need to be modified regarding the calculation of
\g{\mm G_n}. In fact, all nonlinear element formulations compatible
with the HBM, including time-discrete Alternating-Frequency-Time
schemes \zo{guil1998,naci2003} and event-driven frequency-domain
schemes \zo{petr2003,krac2013b}, are applicable. In the authors'
opinion, however, the major advantage of this strategy lies in the
fact that the nonlinear terms are consistent with the steady state
in which the nonlinear forces are also periodic. This provides high
accuracy for the subsequent synthesis procedure, see \sref{nlrom}.\\
It has to be investigated whether this strategy induces inaccuracies
in case of strongly damped systems. In particular, it should be
evaluated whether this strategy for the evaluation of the nonlinear
terms leads to a degradation of accuracy in the prediction of
transient system behavior.

\subsection{Condensation of the eigenvalue problem\label{sec:ana_condensation}}
In many cases the nonlinear forces \g{\mm g(\mm u,\mm{\dot u})} in
\eref{eqm_autonomous} and its Jacobian are highly sparse. By
exploiting this sparsity, the computational effort for the modal
analysis can be significantly reduced, particularly in case of
localized nonlinearities. This strategy has already been followed by
several researchers for forced response analysis using the HBM, see
\eg \zo{chen1998,yang1999a,petr2003,naci2003}. In this study, it is
applied for the first time to the modal analysis of autonomous
systems. We propose to consider the spectral decomposition of the
linearized system for this task, which is particularly beneficial in
this case, as detailed later.\\
The system of equations \erefo{eqm_autonomous} can be partitioned
with respect to \g{\nnl} nonlinear \myquote{N} and
\g{\nll=\ndim-\nnl} linear \myquote{L} terms,
\e{\mm g = \vector{\mm g^{\rm N}\left(\mm u^{\rm N}\right)\\ \mm 0}
\quad\text{with}\quad \mm u = \vector{\mm u^{\rm N}\\ \mm u^{\rm
L}}\fp}{u_partition}
Reformulating \eref{complex_evp} accordingly and premultiplying with
the dynamic compliance matrix for each harmonic $n$, \g{\mm H_n=\mm
S_n^{-1}} yields
\e{\vector{\mm U^{\rm N}_n\\ \mm U^{\rm L}_n} + \matrix{cc}{\mm H^{\rm{NN}}_n & \mm H^{\rm{NL}}_n\\
\mm H^{\rm{LN}}_n & \mm H^{\rm{LL}}_n}\vector{\mm G^{\rm N}_n\\ \mm
0} = \vector{\mm 0\\ \mm 0}\,,\quad n = 0,\cdots,\nh\fp}
{complex_evp_partition}
As the nonlinear forces only depend on the nonlinear unknowns, it is
sufficient to solve only the nonlinear part iteratively,
\e{\mm U^{\rm N}_n + \mm H^{\rm{NN}}_n\mm G^{\rm N}_n\left(\mm
U_0^{\rm N}\,,\,\cdots\,\mm U_{\nh}^{\rm N}\right)=\mm 0\,,\quad n =
0,\cdots,\nh\fp}{complex_evp_condensed}
The dimension of the system of equations in \eref{complex_evp},
which is proportional to the number of coordinates of the full
system, \g{\ndim}, can therefore be reduced to the dimension of
\eref{complex_evp_condensed}, which is proportional to the number of
coordinates associated to nonlinear elements, \g{\nnl\ll\ndim}. If
required, the remaining DOFs can be easily recovered using \g{\mm
U^{\rm L}_n=-\mm H^{\rm{LN}}_n\mm G^{\rm N}_n}. This expansion is
required to evaluate the normalization conditions given by
\erefss{amplitude_normalization} or \erefo{energy_normalization}.\\
It should be noted that the factorization of the dynamic stiffness
matrix has to be computed in each iteration of the nonlinear solver
as it depends on the unknown eigenvalue \g{\lambda}. This can in
general diminish the advantage of this condensation. Owing to the
monomial form of \g{\mm S} in \eref{dyn_stiff}, however, the
inversion can be accomplished very efficiently by using the spectral
decomposition of the structural matrices
\e{\mms\phi_k^{\rm H}\mm M\mms\phi_k = 1\,,\,\, \mms\phi_k^{\rm
H}\mm K\mms\phi_k = \omega_k^2\,,\quad k = 1,\cdots,\ndim\fp}{evlin}
The expensive matrix inversion can then be restated as a simple
matrix product and the trivial inversion of a diagonal matrix,
\e{\mm H_n(\lambda) =
\suml{k=1}{\ndim}{\frac{\mms\phi_k\mms\phi_k^{\rm H}}
{\omega_k^2+(n\lambda)^2}}\,,\quad n = 0,\cdots,\nh\fp}{hinv}
As mentioned in \sref{cnma}, the stiffness matrix \g{\mm K} contains
the linear part of the system. The linearized modal basis in
\eref{evlin} only has to be computed once and for all. It can then
be used as a starting guess for the NMA and incorporated in the
efficient condensation technique proposed in \eref{hinv}. The
availability of the linearized modes has the additional advantage
that the NMA only has to be carried out in the actually nonlinear
regime. In fact, it is important to notice that the expression used
in \eref{hinv} cannot be used in the linear regime since then the
denominator corresponding to \g{n=1} vanishes.

\subsection{Extension to systems of general structure\label{sec:ana_general}}
The proposed method can easily be applied to more general
second-order systems, \ie with linear symmetric and skew-symmetric
velocity- and displacement-dependent terms in \eref{eqm_autonomous}.
The dynamic stiffness matrix in \eref{dyn_stiff} then has to be
augmented accordingly. The spectral decomposition in \eref{evlin} is
obtained by solving a general, quadratic eigenvalue problem in this
case. This decomposition can still be used to assemble the inverse
dynamic stiffness matrix in \eref{hinv}. This decomposition is
derived in \aref{hgen}. It should, however, be noted that linear
damping is proposed to be accounted for in the synthesis rather than
in the modal analysis, so that it can be varied without
re-computation of modal properties, see \sref{nlrom}.

\subsection{Numerical aspects\label{sec:ana_num}}
The complex eigenvalue problem in \eref{complex_evp} combined with
appropriate normalization conditions represents a system of
nonlinear algebraic equations, which has to be solved simultaneously
within a specified energy range. The result of this solution process
are the energy-dependent nonlinear modes, corresponding
eigenfrequencies \g{\omega_0} and modal damping ratios \g{D}. Of
course, the energy range has to cover the range in which the
response of the system is of interest. The synthesis procedure
proposed in \sref{nlrom} is therefore restricted to the energy range
for which the modal properties have been computed.\\
The resulting system of equations was solved using a Newton-Raphson
method. The eigensolution of the linearized system was taken as an
initial guess for a small energy level. In contrast to \zo{laxa2009}
the solution was continued using a predictor-corrector continuation
scheme, see \eg \zo{seyd1994}. This continuation was necessary to
compute the complex, multi-valued relationship between the nonlinear
modal properties and the energy often reported in this context
\zo{lee2005,kers2009}. It should be noticed that more elaborate
bifurcation and stability analysis methods would represent ideal
complements to the framework addressed in this study.\\
The computational efficiency of the solution process was greatly
enhanced by providing analytically calculated derivatives of
\erefs{complex_evp_condensed}{energy_normalization} with respect to
the unknown harmonic components of the eigenvector and the complex
eigenvalue. Analytical derivatives were obtained from manual
symbolic differentiation for each type of nonlinearity, as described
\eg in \zo{siew2009a}. Automatic differentiation could generally
also be used for this task, see
\zo{krac2013b}.\\
The continuation of the solution generally has to be performed for a
large range of the modal amplitude. It is typically not a priori
known in what energy ranges the most relevant regimes of the system
are. In this study, a logarithmic scaling of the modal amplitude
resulted in great computational efficiency in this regard. A linear
scaling of the various unknowns (displacement, frequency, damping,
energy) was applied in order to obtain approximately matching orders
of magnitude, which can have a crucial influence on the convergence
behavior of the nonlinear solver.

\section{Nonlinear modal synthesis\label{sec:nlrom}}
In a linear system, it is possible to formulate the general response
as a synthesis of all solutions to the eigenproblem. This synthesis
is very efficient owing to the superposition principle and the
orthogonality conditions between the modes. In a nonlinear system,
however, these relationships do not hold anymore so that further
assumptions and restrictions have to be accepted for any type of
synthesis procedure. In this study, we restrict the synthesis to the
periodic steady-state vibrations of harmonically forced and
self-excited systems and assume the absence of internal resonances.
The latter aspect inherently excludes systems which already have
multiple
eigenvalues in the linear case.\\
The equation of motion now takes the form
\e{\mm M \mm{\ddot u} + \mm C\mm{\dot u} + \mm K\mm u + \mm g(\mm u,
\mm{\dot u}) = \Re{\mm f_1\ee^{\ii\Omega t}}\fp}{eqm}
Compared to the autonomous case in \eref{eqm_autonomous}, the
equation is augmented by a real symmetric viscous damping matrix
\g{\mm C} and a forcing term of complex amplitude \g{\mm f_1} and
frequency \g{\Omega}. This forcing term vanishes in the self-excited
setting.\\
In order to solve \eref{eqm} by nonlinear modal synthesis, we apply
the single-nonlinear-resonant-mode theory \zo{szem1979}. This theory
is based on the observation that in the absence of nonlinear modal
interactions, the energy is basically concentrated in a single mode
'$j$'. This mode dominates the system response and is treated as
nonlinear. Owing to their low energy level, the remaining modes are
approximated by their linearized counterpart \g{\mms\phi_k} in
accordance with \eref{evlin}. With this assumption, the system
response can be formulated as
\e{\mm u(t) \approx \Re{
\underbrace{q_j\suml{n=0}{\nh}{\mms\psi_n(\left| q_j \right|
)\ee^{\ii n\Omega t}}}_{\text{nonlinear mode $j$}} +
\underbrace{\suml{k\ne j}{\ndim}{q_k\mms\phi_k\ee^{\ii\Omega
t}}}_{\text{linearized modes}} } \fp}{snrm}
%
Herein, the fundamental frequency \g{\Omega} is either the
excitation frequency of the harmonic excitation or the frequency of
the self-excited vibration, \g{\Omega=\omega_j}. The mode number $j$
has to be selected so that either the excitation is around the $1:1$
external resonance of this mode or the self-excitation leads to a
limit-cycle-oscillation (LCO) in this mode. A new modal amplitude
\g{q_j} and eigenvector with complex amplitudes \g{\mms\psi_n} have
been introduced in \eref{snrm}. The relationship to the variables
obtained from the modal analysis is defined as follows:
\e{\mms\psi_1^{\rm H}\mm M\mms\psi_1=1 \quad \Rightarrow \quad
q_j\mms\psi_n = \mm U_n\,,\,\, n =
0,\cdots,\nh\fp}{mass_normalization}
For clarity, the formal dependence of \g{\omega_j,D,\psi_j} on
\g{\left| q_j \right|} is not denoted here and in the following.\\
Calculation of the modal amplitudes in \eref{snrm} is carried out by
once again employing single-nonlinear-resonant-mode theory: The
contributions of the linearized modes are calculated in the
traditional manner, \ie by projecting the linearization of
\eref{eqm} onto the linear part of the modal basis. The modal
amplitude of the nonlinear mode is determined independently of the
linear modes by projecting the equation of motion formally onto each
harmonic \g{\mms\psi_n\ee^{\ii n\Omega t}} of the nonlinear
eigenvector. The fundamental harmonic \g{n=1} yields
\e{\left[-\Omega^2 + \ii\Omega\mms\psi_1^{\rm H}\mm
C\mms\psi_1+\omega_j^2+2D_j\omega_j~\ii\Omega\right]q_j =
\mms\psi_1^{\rm H}\mm f_1\fp}{eqm_projected}
The last two terms in the brackets correspond to the projection of
the stiffness matrix and the nonlinear terms onto the $j$-th mode
and are readily available from the modal analysis in \sref{cnma}.\\
It can be easily verified that the resulting non-fundamental
harmonic equations, \ie with \g{n\neq 1} essentially give
\g{\mms\psi_n^{\mathrm H}\left(\mm S_n\mm U_n + \mm G_n\right) = 0},
which is inherently fulfilled in accordance with the eigenproblem
given by \eref{complex_evp}. This applies exactly in case of
conservative systems in resonance. For dissipative systems and/or
systems driven not precisely at resonance, this approach represents
an approximation: Firstly, the dynamic stiffness matrix \g{\mm S_n}
used in \eref{complex_evp} is not identical to the one associated
with \eref{eqm} in this case. Particularly, it is evaluated at
\g{\Omega} rather than at \g{\omega_j}. Similarly, the nonlinear
forces are evaluated at the excitation frequency rather than the
eigenfrequency. By assuming that these effects can be neglected in
the vicinity of the resonance, the multi-harmonic response in
\eref{snrm} can be synthesized without the need to solve any
nonlinear equations in addition to \eref{eqm_projected}. This is in
contrast to \zo{laxa2009}, where only the fundamental harmonic of
the nonlinear mode was considered instead of the multi-harmonic,
multi-modal response in \eref{snrm}. It is conceivable that the
accuracy of the approach could generally be increased by overcoming
these simplifications, but this would most certainly result in
additional computational effort.\\
Note that \g{\psi_1,D_j,\omega_j} in \eref{eqm_projected} depend on
the modulus \g{\left|q_j\right|} of the complex modal amplitude.
Therefore, \eref{eqm_projected} must be solved iteratively. Owing to
the numerical character of the modal analysis, the modal properties
will only be available at discrete amplitude values. Hence, a
one-dimensional interpolation scheme was used in order to apply the
continuous formulation in \eref{eqm_projected} to the numerical
results of the modal analysis described in \sref{cnma}. It was found
that both piecewise linear and piecewise cubic interpolation schemes
performed well in this study.

\paragraph{Frequency response function synthesis}
The frequency response function (FRF) can be obtained by solving
\eref{eqm_projected} and employing a continuation on \g{\Omega}.
Once the modulus of the modal amplitude \g{\left|q_j\right|} has
been computed, the phase of the nonlinear modal amplitude is
obtained by evaluating \eref{eqm_projected}. The modal amplitudes of
the remaining, linearized modes are also calculated in the
traditional manner.

\paragraph{Backbone curve synthesis}
Often not the whole frequency response function is relevant but only
the backbone curve, \ie the direct relationship between resonance
amplitude, resonance frequency and a system parameter. This backbone
curve can be obtained from \eref{eqm_projected} by setting
\g{\Omega=\omega_j} and employing a numerical continuation on the
desired system parameter. Similar to the FRF case, the phase and
contribution of linearized modes have to be determined to synthesize
full response in \eref{snrm}.

\paragraph{Calculation of self-excited limit cycles}
In case of self-excited vibrations, the right hand side of
\eref{eqm_projected} is zero. It directly follows that
\g{\Omega=\omega_j} and the equation simplifies to
\e{ \mms\psi_1^{\rm H}\mm C\mms\psi_1 + 2D_j\omega_j = 0\fp}{lco}
\eref{lco} governs the limit cycle oscillation amplitude
\g{\left|q_j\right|} of the nonlinear mode. In the self-excited
case, the phase is not relevant so that only the modulus needs to be
determined. The assumption of a single nonlinear resonant mode
implies that the linearized modes must be damped away. Otherwise,
they would grow unbounded and/or lead to nonlinear modal
interactions, which have been excluded from the synthesis in this
study. The full response is thus recovered by the first term in
\eref{snrm}.

\subsection{Advantages of the proposed nonlinear modal synthesis\label{sec:nlrom_adv}}
A major benefit of the proposed methodology is that only a scalar
nonlinear equation has to be solved for the nonlinear modal
amplitude, independent of the number of DOFs originally contained in
the system, the type and distribution of nonlinear sources or the
number of harmonics in the modal analysis. Hence, the computational
effort for the synthesis procedure is almost negligible compared to
alternative methods for the computation of the steady state dynamics
such as HBM or
direct time integration in conjunction with shooting.\\
It should also be noticed that the approximated response in
\eref{snrm} is still multi-harmonic and multi-modal in contrast to
the approach proposed in \zo{laxa2009}. The multi-harmonic character
is an important characteristic of the time evolution of the
nonlinear forces and the resulting response. The method that we
propose here, models this important
feature of the solution.\\
Moreover, the prediction of steady-state vibrations can be carried
out in a large range of parameters of the system in \eref{eqm},
without the need for expensive re-computation of the modal basis.
This provides enormous computational savings, in particular for
exhaustive parameter variations that are typically required \eg for
uncertainty analysis and optimization.

\subsubsection{Variation of the excitation}
Since the nonlinear modes were calculated independent of the
excitation force distribution, any force field can be applied to the
system. It is only important that the modes have been computed
within the energy range, in which the system is driven by the
excitation.\\
This is particularly interesting for applications where the load
collective - in terms of spatial distribution and amplitude and phase - is uncertain.
Exhaustive parametric studies on the forced response can then be
conducted at virtually no extra computational effort.

\subsubsection{Damping\label{sec:damping}}
Although viscous damping was considered initially, other common
types of damping can be studied by simply replacing or augmenting
the term \g{\ii\Omega\mms\psi_1^{\rm H}\mm C\mms\psi_1} in
\eref{eqm_projected}. For hysteretic damping \g{\mm D_{\rm{hyst}}},
this term results in \g{\ii\mms\psi_1^{\rm H}\mm
D_{\rm{hyst}}\mms\psi_1}. In case of modal damping \g{\eta_j}, the
corresponding term is simply given by \g{\ii\Omega\eta_j}. Similar
to the excitation, damping is often described by uncertain
parameters which have to be varied in the design process of
industrial applications.

\subsubsection{Similar parameter sets}
In some cases, even parameters of the nonlinear term \g{\mm g(\mm u,
\mm{\dot u})} can be varied without re-computation of the nonlinear
modal properties. In this study, this is the case for the considered
systems in \ssref{beam}-\ssref{shroud} where the only nonlinear
terms stem from the preloaded piece-wise linear contact constraints,
\ie the unilateral elastic and the elastic Coulomb contact
constraints. It is postulated that the nonlinear dynamic compliance
is only dependent on the ratio between preload \g{f_{\rm{pre}}} and
excitation level \g{\|\mm f_1\|}. Hence, the response at a different
preload \g{\tilde f_{\rm{pre}}} can be calculated by
\e{q_j\left(\tilde f_{\rm{pre}},\mm f_1 \right) = \frac{\tilde
f_{\rm{pre}}}{f_{\rm{pre}}}
q_j\left(f_{\rm{pre}},\frac{f_{\rm{pre}}}{\tilde f_{\rm{pre}}}\mm
f_1\right)\fp}{preload_excitation_level}
A strict mathematical proof is beyond the scope of this study.
Instead, the performance and accuracy of this hypothesis will be
demonstrated in \sref{numerical_examples}. It should be noticed that
this hypothesis is not required for the general methodology proposed
in this study, but it only provides a beneficial feature for the
specific contact nonlinearities used in \ssref{beam}-\ssref{shroud}.

\section{Numerical examples\label{sec:numerical_examples}}
In this section, the nomenclature regarding the nonlinear modes and
their interaction is similar to the one used by \z{Kerschen
\etal}{kers2009}. A $N:M$ resonance of the $J$th mode is denoted as
\g{\rm{S}N:M,m=J}. If not otherwise specified, the eigenfrequency,
\ie the frequency of the $1:1$ internal resonance is illustrated in
the frequency energy plot for each mode.\\
In all figures, normalized excitation frequencies \g{\Omega^*} and
eigenfrequencies \g{\omega_0^*} are illustrated. The scaling factor
is the linearized eigenfrequency of the first mode, unless otherwise
specified. Furthermore, response amplitudes \g{a^*} have been
non-dimensionalized by scaling with the corresponding linear case.\\
For investigation of the numerical examples, a predictor-corrector
scheme with tangential predictor and pseudo-arc-length
parametrization was implemented and used by the authors. For the
numerical evaluation of the Fourier coefficients of the nonlinear
forces, the well-known time-discrete Alternating-Frequency-Time
scheme was employed, see \eg \zo{came1989}.
\subsection{Modal analysis of a system with cubic spring}
\myf[htb]{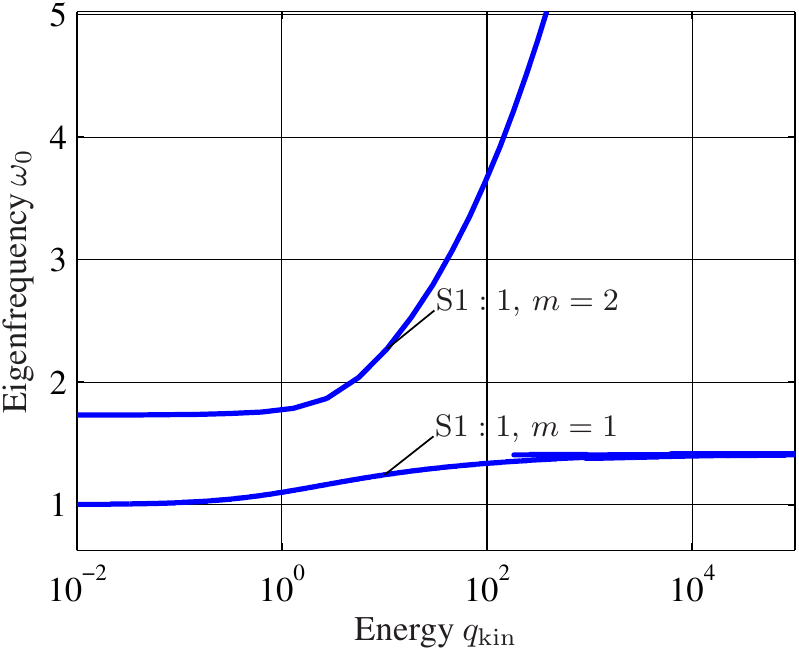}{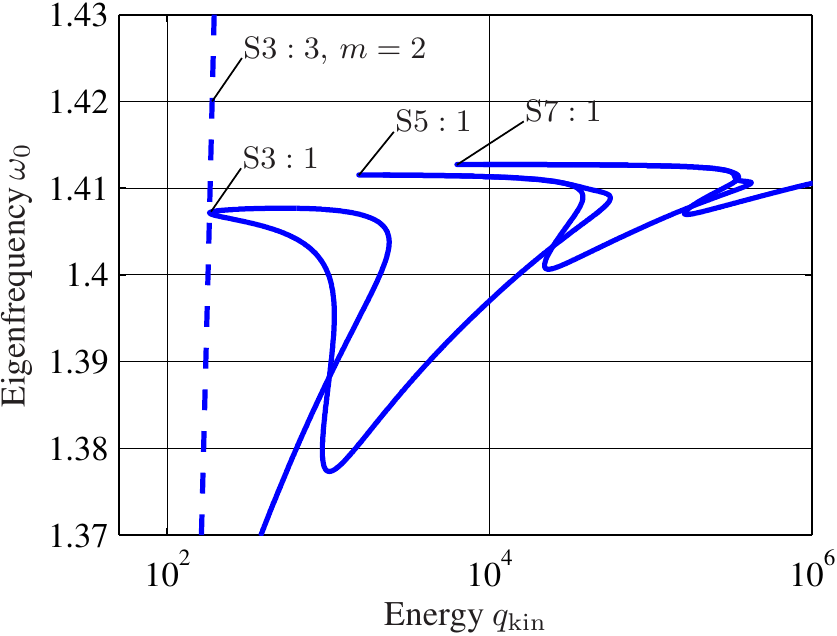}{}{}{.45}{.45}{Frequency-Energy-Plot of
a 2-DOF system with cubic spring (~(a) overview, (b) zoom
on the first internal resonances~)}
As a first validation, the modal analysis technique is first applied
to a 2-DOF system with cubic spring thoroughly studied in
\zo{kers2009}. The equations of motion read
\ea{\nonumber \ddot x_1 + 2 x_1 - x_2 + 0.5 x_1^3 &=& 0\\
\ddot x_2 - x_1 + 2 x_2 &=& 0\fp}{eqm_2dof_cubic}
The frequency-energy plot (FEP) of this system is depicted in
\frefs{fig01a}-
\frefo{fig01b}. The system basically
has two eigenfrequencies corresponding to the in-phase \g{m=1} and
out-of-phase \g{m=2} mode. The stiffening behavior of the cubic
spring becomes apparent.\\
Internal resonances occur when the eigenfrequencies of the two modes
are commensurable at the same energy level. This happens despite the
fact that the linearized eigenfrequencies are not commensurable,
because of the general frequency-energy dependence of each mode. The
3:1, 5:1 and 7:1 internal resonances between the modes are
illustrated in \fref{fig01b}. At the
tip of the \myquote{tongues}, the first eigenfrequency is precisely
\g{\frac13,\,\frac15,\,\frac17} of the second eigenfrequency.\\
In contrast to the time integration scheme in conjunction with a
Shooting procedure proposed in \zo{kers2009}, the modal analysis
technique presented in \sref{cnma} was employed. The results of both
methods are fully equivalent. Several harmonics have to be retained
in the multi-harmonic expansion in order to accurately resolve the
internal resonances. A more detailed bifurcation and stability
analysis of the nonlinear modal interactions was considered beyond
the scope of this paper. The interested reader will find a detailed
analysis of this system is presented in \zo{kers2009}.

\subsection{Analysis of a clamped beam\label{sec:beam}}
\fss[htb]{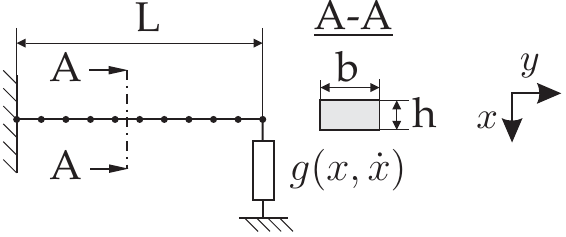}{Clamped beam with nonlinear element}{1.0}
The clamped beam as depicted in \fref{fig02} was investigated.
Its specifications are: Size \gehg{L}{200}{mm}, \gehg{b}{40}{mm},
\gehg{h}{3}{mm}, Young's modulus \gehg{E}{210,000}{MPa} and density
{\x{${\rho}= 7800 ~ {\frac{\rm{kg}}{\rm{m}^3}}$}}. The beam was
discretized by ten Euler-Bernoulli beam elements, the displacement
was constrained to the transverse ($x$) direction. The beam is
connected to nonlinear force elements at its free end. Different
types of nonlinearity \g{g(x,\dot x)} will be investigated with
respect to their effect on the overall vibration behavior of the
system. In the following, the response amplitude \g{a} is defined as
the maximum value of the zero-mean tip displacement.

\subsubsection{Unilateral spring nonlinearity}
\myf[htb]{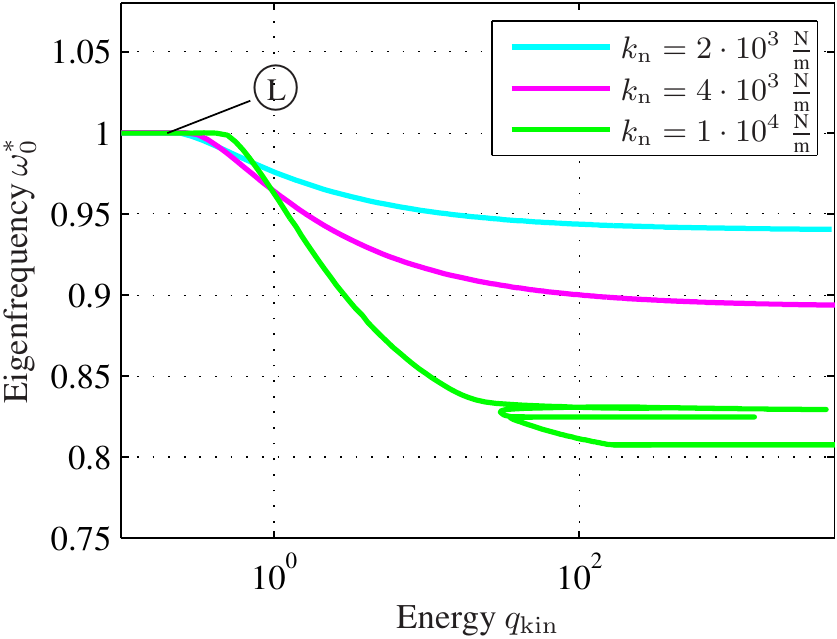}{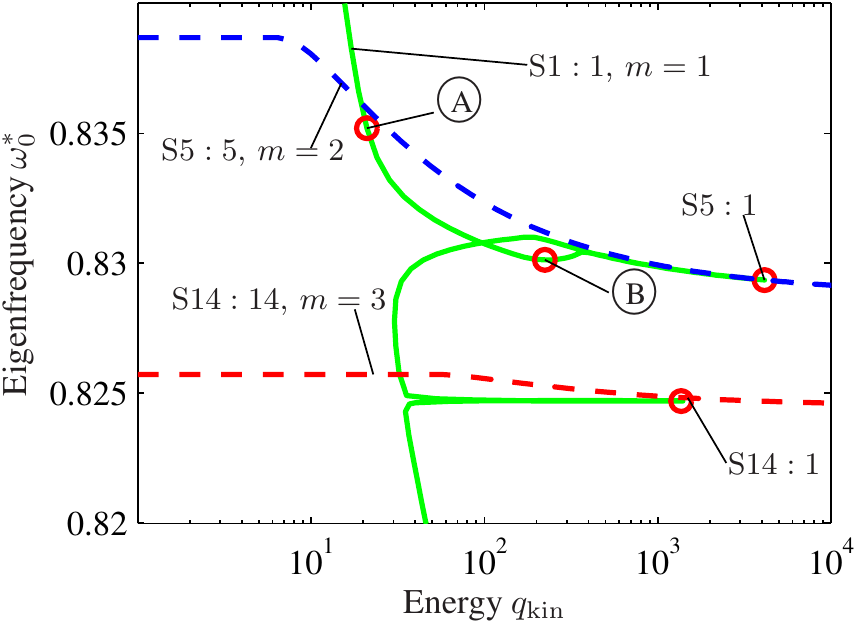}{}{}{.45}{.45}{Frequency-Energy-Plot of a
clamped beam with unilateral preloaded spring (~(a) overview
for different normal stiffness values, (b) zoom on the first internal
resonances for {\x{${\kn}= 1\cdot10^4 ~
{\frac{\rm{N}}{\rm{m}}}$}}~)}
First, a unilateral spring with stiffness \g{\kn} will be
considered,
\e{g(x,\dot x) = \kn\left(x +
a_0\right)_+\fp}{fnl_unilateral_spring}
The value in the parenthesis is only considered if it is greater
than zero. Note that the spring is preloaded by a compression of
length \g{a_0}.\\
The frequency-energy-plot of the first bending mode is depicted in
\frefs{fig03a}-\frefo{fig03b}.
For low energies, the system exhibits linear behavior, \ie the
eigenfrequency remains constant. Once the vibration amplitude is
large enough, there is partial lift-off during the vibration
cycle, causing the apparent softening behavior.\\
In \fref{fig03a} the spring
stiffness value \g{\kn} is varied. The relative frequency shift
increases with the spring stiffness. For large stiffness values, \eg
\g{\kn=10^4}, the system exhibits internal resonances.\\
The internal resonances are indicated in
\fref{fig03b}. In
order to illustrate the modal interactions, the corresponding higher
modes \g{m=2,3} are also depicted. They are denoted $N:N$
resonances, where $N$ is the multiple corresponding to the $N:1$
resonance with the first mode. The following two resonance
coincidences are presented: A $5:1$ resonance with the second mode
and a $14:1$ resonance with the third mode.
\mythreef[t!]{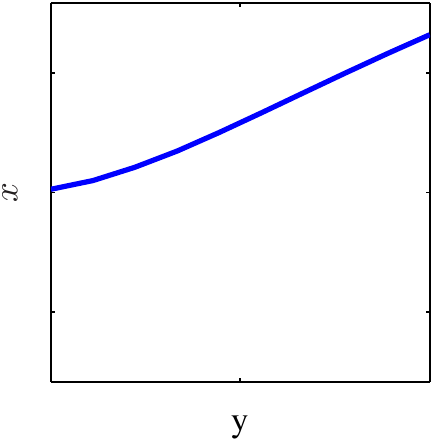}{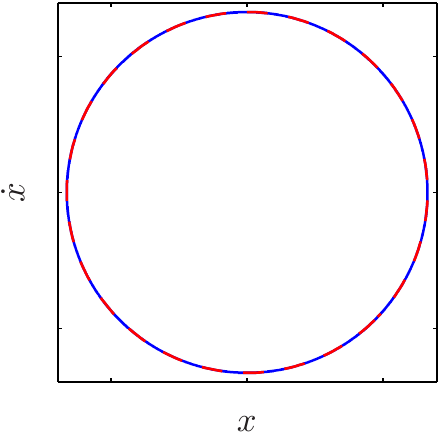}{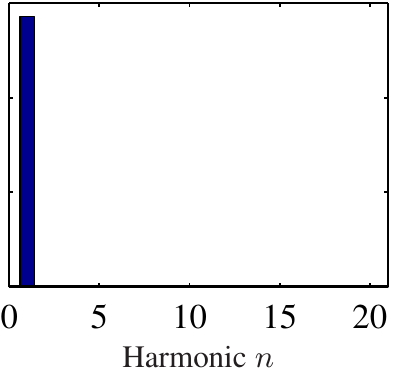}
{}{}{}{Linear case
\kreis{L} (~(a) mode shape, (b) phase portrait, (c) frequency content~)}
\mythreef[h!]{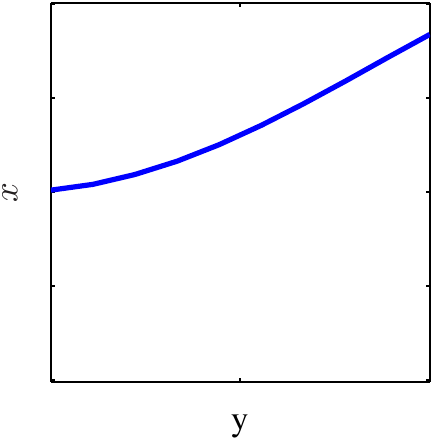}{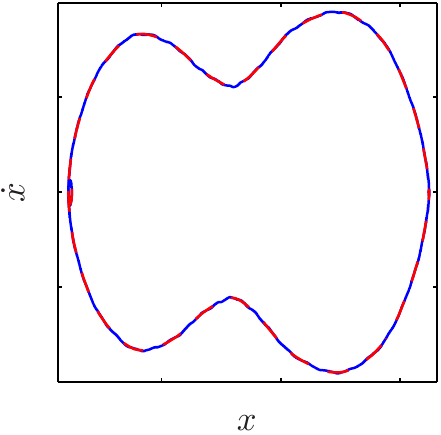}{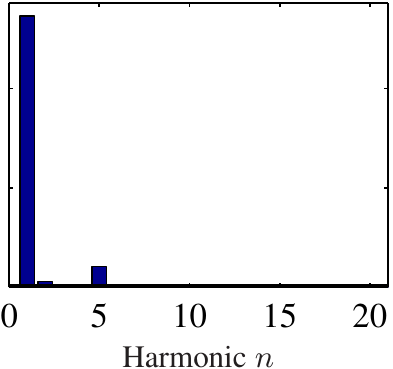}
{}{}{}{Case \kreis{A} (~(a) mode shape, (b) phase portrait, (c) frequency content~)}
\mythreef[h!]{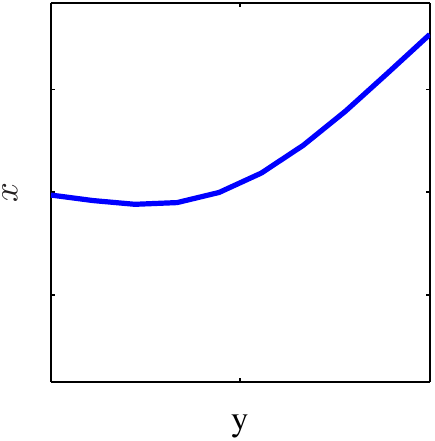}{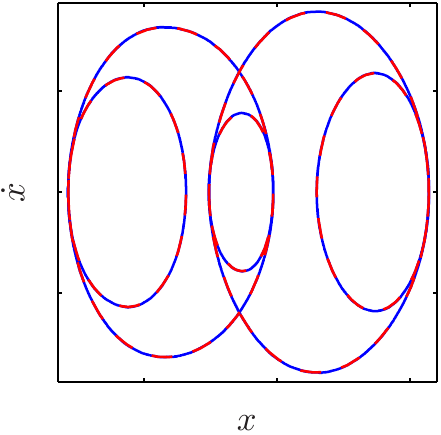}{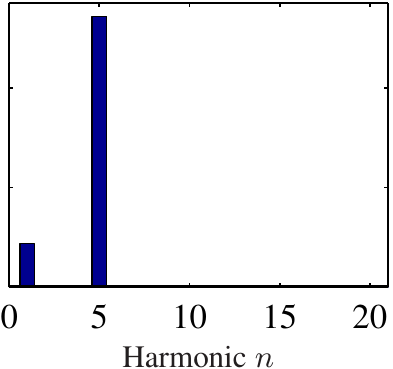}
{}{}{}{Case \kreis{B} (~(a) mode shape, (b) phase portrait, (c) frequency content~)}
\mythreef[t!]{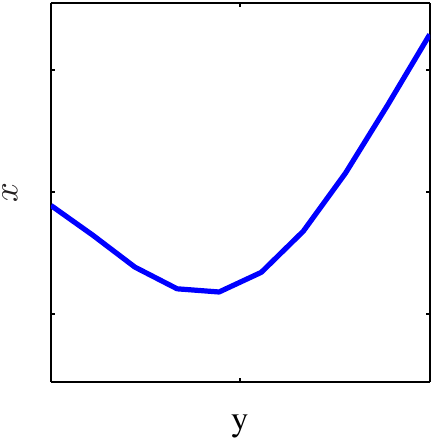}{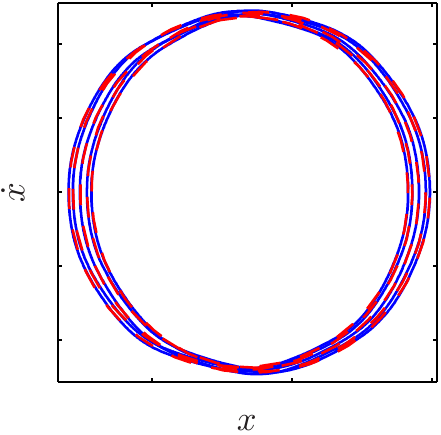}{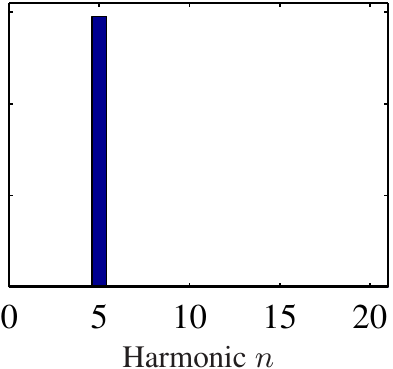}
{}{}{}{Internal resonance
S5:1 (~(a) mode shape, (b) phase portrait, (c) frequency content~)}
\mythreef[h!]{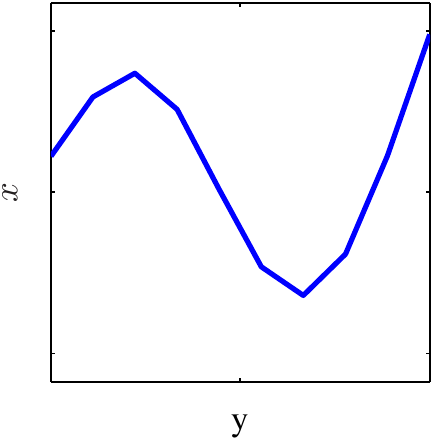}{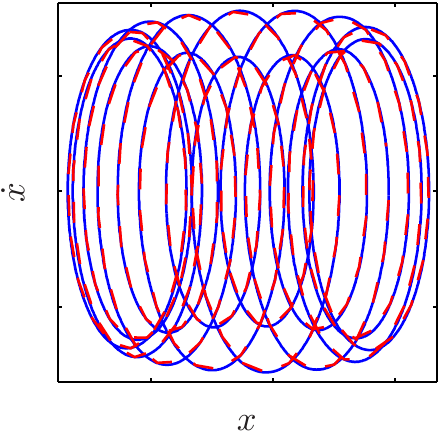}{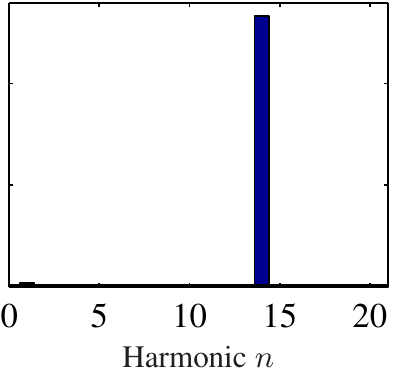}
{}{}{}{Internal resonance
S14:1 (~(a) mode shape, (b) phase portrait, (c) frequency content~)}
\\
In order to provide a better understanding of the underlying
dynamics, the nonlinear modes for the points indicated in
\frefs{fig03a}-\frefo{fig03b}
are depicted in
\frefs{fig04a}-\frefo{fig08b}.
Mode shape, phase portraits for the beam tip deflection and
frequency content of the kinetic energy are illustrated. In vicinity
of the internal resonance, \ie at the tip of the tongues, an
according higher harmonic content becomes apparent while the
frequency content far from these points is dominated by the
fundamental harmonic. As it can be deduced from the figures, also
the mode shape in the vicinity of these resonances becomes similar
to the interacting mode shape.\\
The complex system behavior is well-resolved by the proposed modal
analysis technique, which can be ascertained by comparing the phase
portrait results of the proposed method (solid) to the
time-step integration results (dashed), \cf
\frefs{fig04b}-\frefo{fig08b}.
\fss[htb]{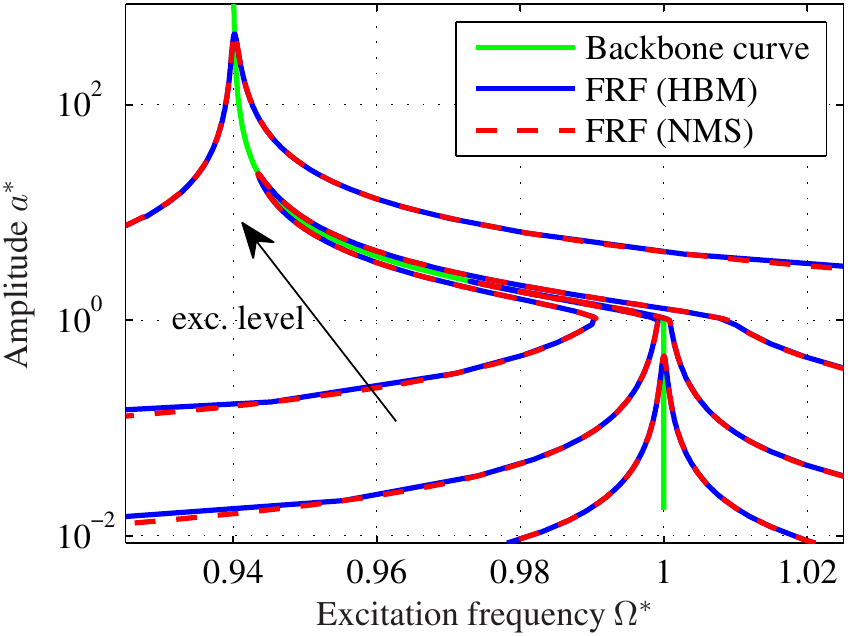}{Forced response of a
clamped beam with unilateral preloaded spring for varying excitation
level}{1.0}
\\
The synthesis method proposed in \sref{nlrom} is restricted to
regimes where internal resonances are absent. Hence, a value of
\g{\kn=2\cdot10^3} was specified for the subsequent investigations.
For this value, the forced response to a discrete harmonic
excitation at the middle of the beam in \fref{fig02} was
calculated. The system was excited in a frequency range around the
eigenfrequency of the first bending mode. A hysteretic damping
\g{\mm D_{\rm{hyst}} =\eta\mm K} with a damping factor of
\g{\eta=0.1\%} was specified, see \ssref{damping}. The results are
depicted in \fref{fig09}. The normalized
amplitude is defined as \g{a^*=\frac{a}{a_0}}.\\
The excitation level has been varied in a wide range. For increasing
excitation level the modal amplitude increases and a softening
effect becomes apparent, in full accordance with the results in
\fref{fig03a}. Note the overhanging
branches resulting in a multi-valued forced response.\\
The forced response has also been computed with the conventional,
multi-term harmonic balance method (HBM) with a harmonic order
\g{\nh=7}. It can be ascertained that the synthesis (NMS) of the
forced response is in very good agreement with HBM results. In
particular in the vicinity of the resonance, the accuracy of the
proposed synthesis is excellent\footnote{The expression 'excellent
accuracy' was used throughout this study to indicate that relative
errors compared to the reference did not exceed $1\%$.}. Moreover
the backbone curve, which has also been directly obtained by the
modal synthesis, matches well with the resonances of the forced
responses.

\subsubsection{Friction nonlinearity}
\myf[htb]{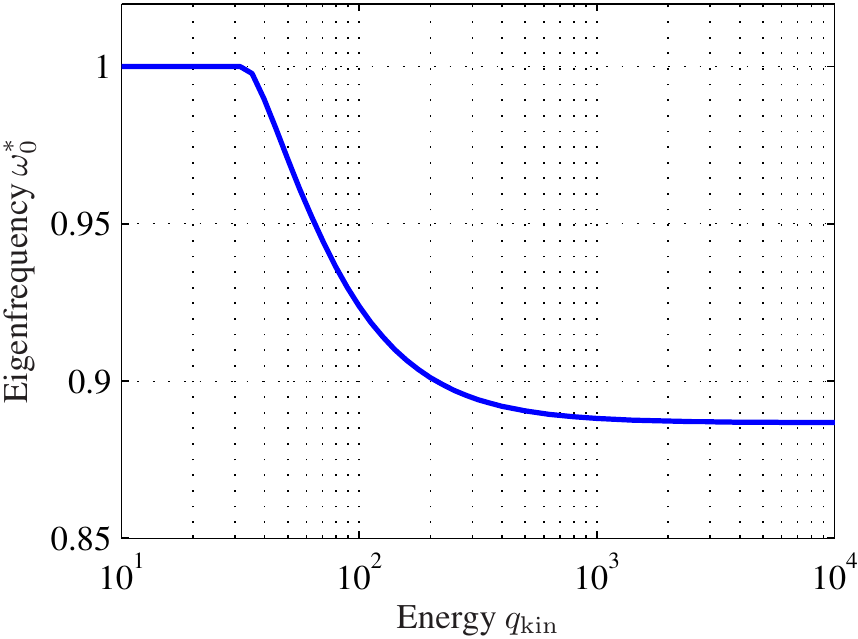}{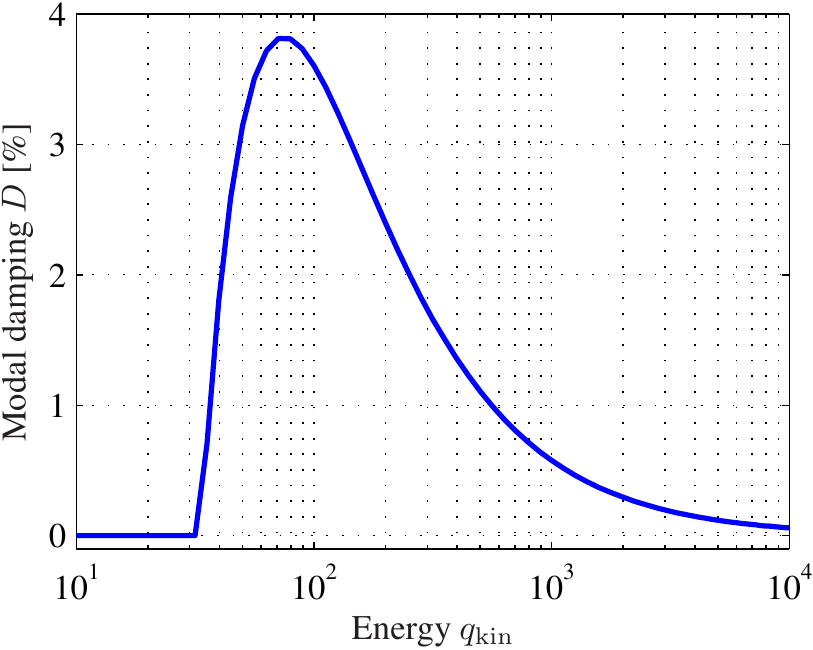}{}{}{.45}{.45}{Modal properties of a clamped beam with friction
nonlinearity (~(a) eigenfrequency, (b) modal
damping~)}
\fss[h!]{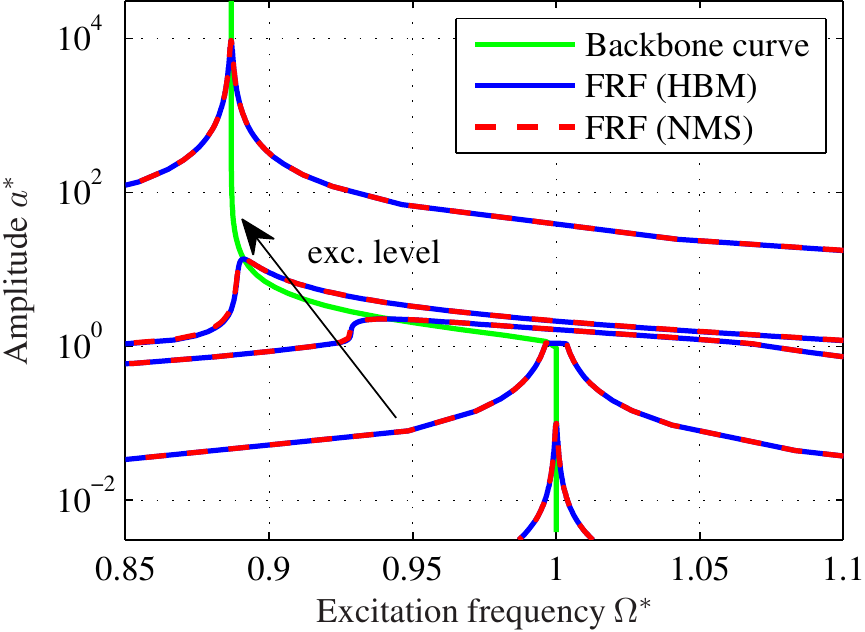}{Forced response of a
clamped beam with friction nonlinearity for varying excitation
level}{1.0}
Next, an elastic Coulomb nonlinearity with a stiffness \g{\kt} and a
limiting friction force \g{\mu N} is considered. The expression for
the nonlinear force is given by the following differential equation
\e{\dot g(x,\dot x) = \begin{cases}0 & \left|g(x,\dot x)\right|=\mu
N\\
\kt\dot x & \left|g(x,\dot x)\right|<\mu N\end{cases}\fp}{fnl_fric}
Several approaches exist for the regularization of this nonlinearity
in time and frequency domain
\zo{yang1998a,guil1998,petr2003}. In this study, the
time-discretized formulation in \zo{guil1998} was used.\\
The modal properties are illustrated in
\frefs{fig10a}-\frefo{fig10b}.
A softening behavior can be ascertained from the FEP. The modal
damping is zero in the fully stuck state. The damping increases to a
maximum in the microslip regime and decreases asymptotically to zero
for large energy levels. The decrease of the modal damping ratio may
seem counter-intuitive at first, but it can be easily made
plausible: For a viscous damping source, the dissipated energy grows
quadratically with amplitude, leading to a constant modal damping
ratio. The energy dissipated in the Coulomb slider essentially
increases only linearly with amplitude, thus leading to a decreasing
modal damping for large amplitudes. The interested reader is
referred to \zo{popp2003a,laxa2009} for further insight in the
qualitative dynamic behavior of structures with friction joints.\\
\paragraph{Forced response synthesis}
In \fref{fig11}, the forced response is
depicted for a varying excitation level. Again, the results of the
proposed synthesis method (NMS) are in excellent agreement with
conventional forced response calculations (HBM).
\paragraph{Calculation of limit cycles}
For the analysis of self-excited vibrations, the viscous damping
matrix \g{\mm C} was defined by inverse modal transformation,
\e{\mm C = \mms\Phi^{-\rm H}\mm{\diag}\lbrace 2
D_k\omega_k\rbrace\mms\Phi^{-1}\fp}{clin_flutter}
Herein, \g{\mms\Phi=\matrix{ccc}{\mms\phi_1 & \cdots &
\mms\phi_\ndim}} and \g{\omega_k} are the modal matrix and the
eigenfrequencies of the linearized system defined in \eref{evlin}.
It should be noted that this damping definition is common for
simplified flutter analyses in turbomachinery applications, see \eg
\zo{petr2012c}. A negative value was specified for the first modal
damping ratio \g{D_1} to obtain self-excited vibrations in the first
mode. The remaining damping ratios were defined as \g{D_k = 1\%,\,
k=2,\cdots,\ndim}.
\fss[t!]{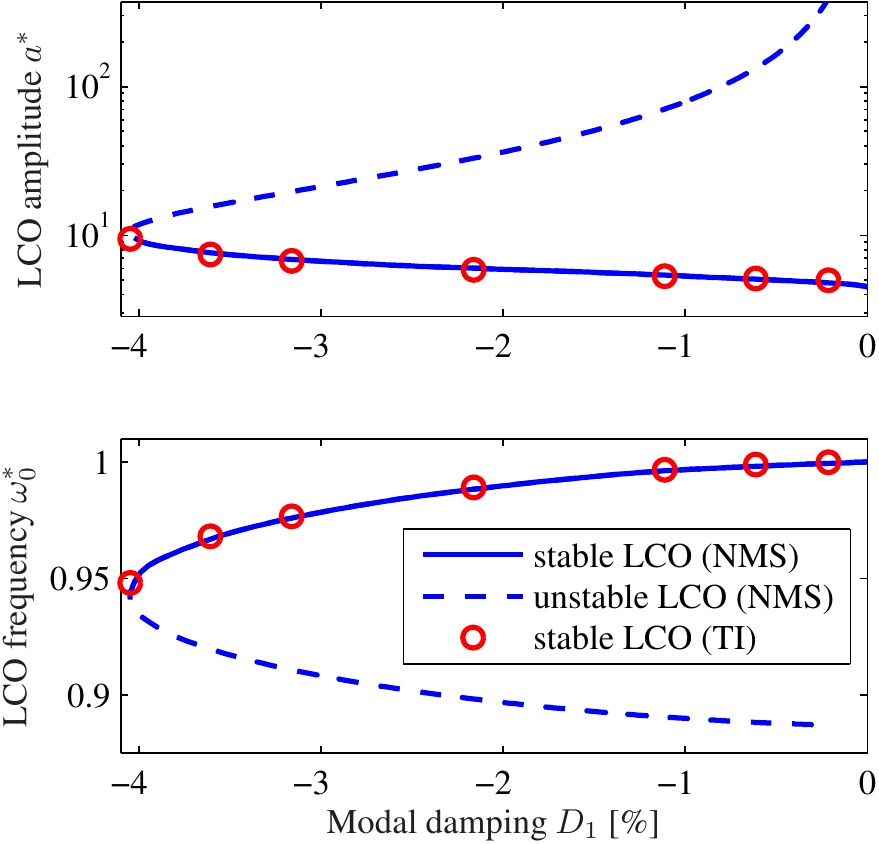}{LCO amplitude and frequency
w.r.t. modal damping ratio \g{D_1}}{1.0}
\\
In \fref{fig12}, the LCO amplitude and frequency
are illustrated with respect to the modal damping \g{D_1}. Stable
and unstable regimes exist. The local stability was determined
simply by considering the sign of the slope of the modal damping at
the limit cycle amplitude: If the modal damping increases with
respect to the amplitude, the limit cycle is stable, otherwise it is
not stable. Whether the system reaches a LCO depends on the initial
energy in the system. For sufficiently low damping values
\g{D_1<-4\%}, an LCO does not exist, \ie the modal amplitude would
grow unbounded in this case. For positive damping values, the LCO
degenerates to the equilibrium point. It should be noted that the
limitation of the vibration amplitude is not only influenced by the
amount of nonlinear modal damping, but it also depends on the mode
shape, see \eref{lco}. The mode shape determines the effective modal
damping \g{\mms\psi_1^{\rm H}\mm C\mms\psi_1}, assuming a constant
damping matrix \g{\mm C}. The stable LCOs computed by the synthesis
(NMS) are in excellent agreement with time-step integration
simulations, see \fref{fig12}.

\subsection{Analysis of a turbine bladed disk with shroud contact\label{sec:shroud}}
\myf[ht!]{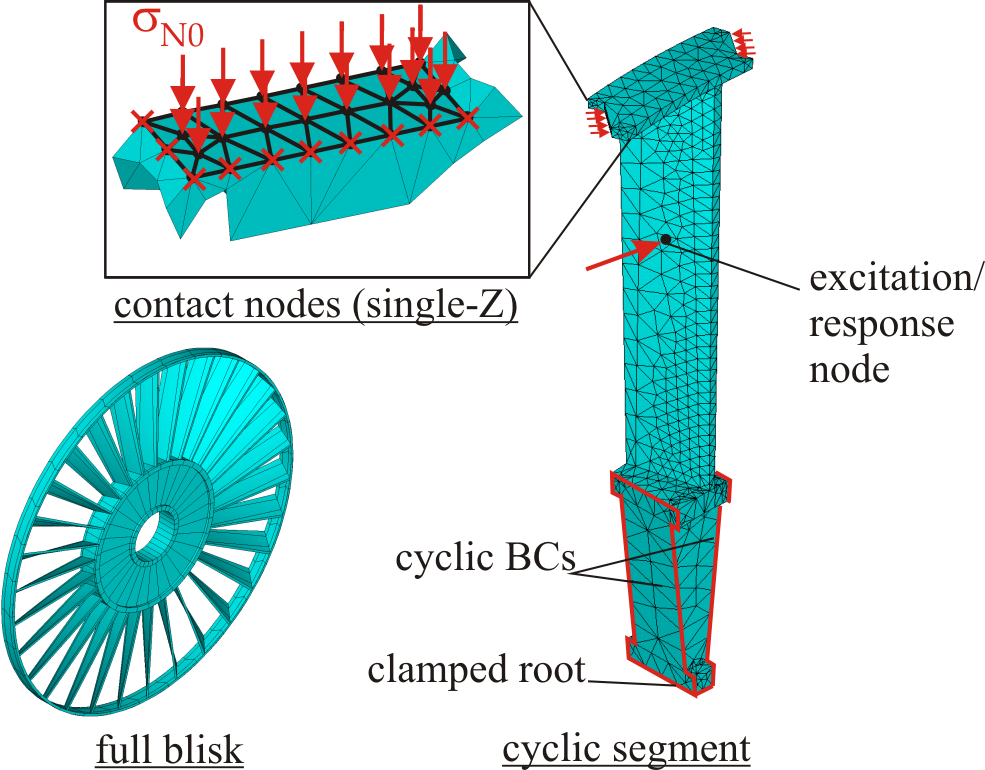}{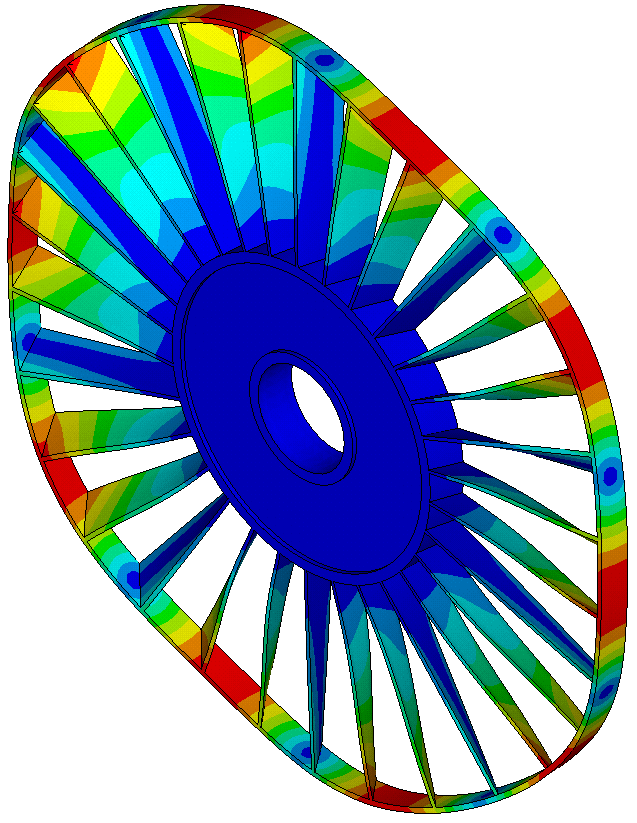}{}{}{.55}{.35}{Finite Element Model of a bladed disk with
shroud contact (~(a) system
definition, (b) investigated mode shape for fixed contact
conditions~)}
In order to demonstrate the applicability of the proposed
methodology to more complex systems, a turbine bladed disk with
shroud contact interfaces is considered. It comprises \g{30} blades
and is considered to be of perfect cyclic symmetry, \cf
\fref{fig13a} and is similar to the ones investigated in
\zo{krac2012a,krac2013c,krac2013g}. Owing to the symmetry, only a single
sector with cyclic boundary conditions was regarded. The finite
element model of the segment consists of \g{25641} DOFs. A cyclic
Craig-Bampton reduced order model \zo{siew2009a} of the sector was
constructed. Only the DOFs involved in the contact formulation were
retained as master coordinates. A number of $50$ linear normal modes
were ascertained to yield convergence of the results obtained in
this study. The second mode family for a spatial harmonic index $5$
was considered in the subsequent investigations, see
\fref{fig13b}. In the forced case, a discrete
traveling-wave-type excitation was specified. The amplitude \g{a} is
defined as the maximum value of the time domain displacement in the
circumferential direction at the response node, \cf
\fref{fig13a}.
\\
Impenetrability and friction constraints were enforced in terms of
unilateral springs in normal direction and elastic Coulomb elements
in the tangential plane \zo{siew2009a,krac2012a}. The
three-dimensional contact model therefore allows for stick, slip and
lift-off phases, and takes into account the influence of the normal
dynamics on the stick-slip transitions. The contact constraints were
imposed in a node-to-node formulation.\\
Contact was defined only in the central part of the Z-shaped shroud
in this study. It is often desirable to design the normal pressure
distribution in such a way that the shroud interfaces are not fully
separated during operation so that no high-energy impacts occur and
the resonance frequency shifts remains small. Hence, it is regarded
as realistic that a portion of the contact area is in permanent
contact. For this example, a portion of the nodes was defined as
bonded, indicated by crosses in \fref{fig13a}. A
homogeneous normal pressure distribution was specified for the
remaining contact area. This specific contact scenario has of course
academic character. In order to improve the model accuracy for a
realistic case study, a nonlinear static analysis should be carried
out, taking into account centrifugal effects on the contact
situation and large deformations for the relevant rotational speed
range \zo{siew2009a}.
\myf[t!]{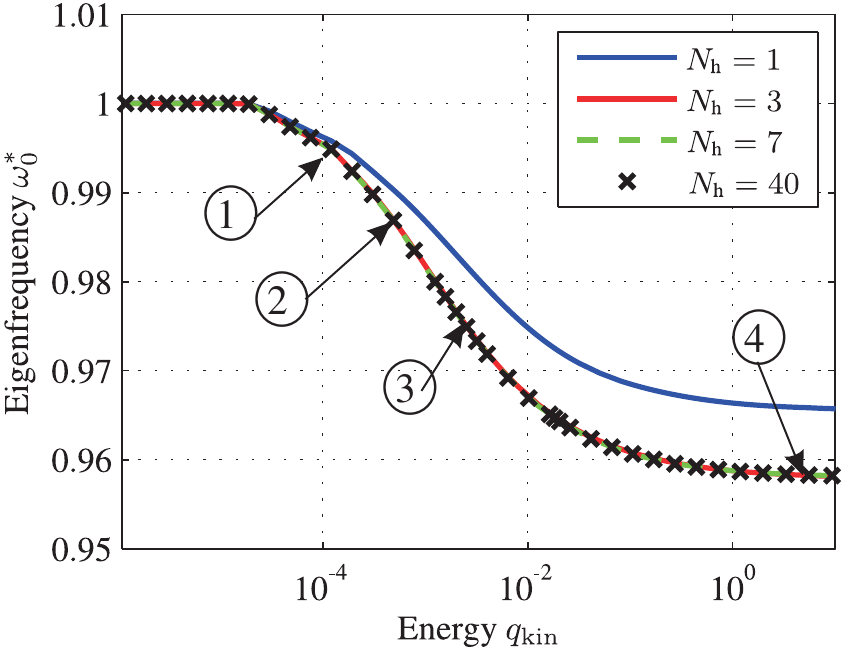}{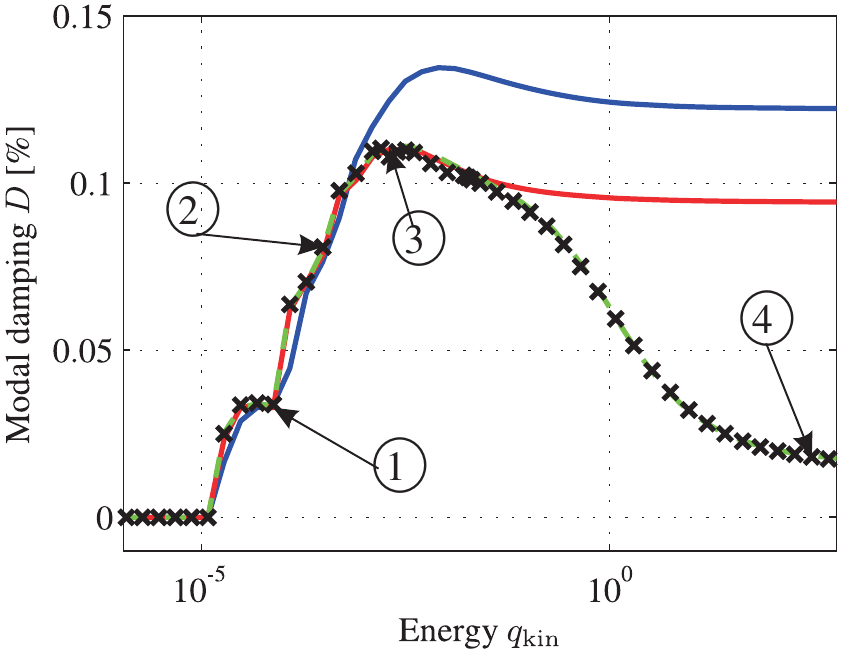}{}{}{.45}{.45}{Modal properties of a bladed disk with shroud
contact (~(a) eigenfrequency, (b) modal
damping~)}
\myf[t!]{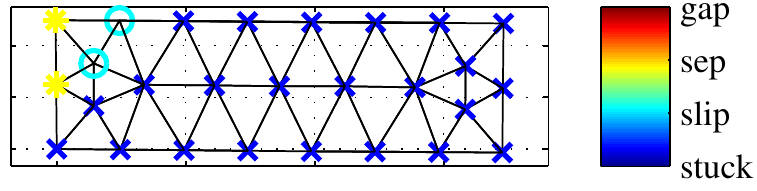}{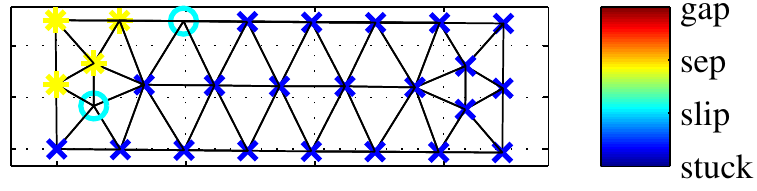}{}{}{.45}{.45}{Contact status for different
modal amplitudes (~(a) point \kreis{1}, (b) point \kreis{2}~)}
\myf[t!]{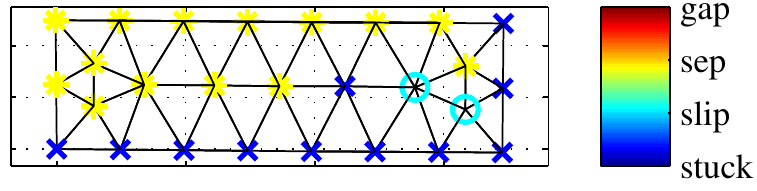}{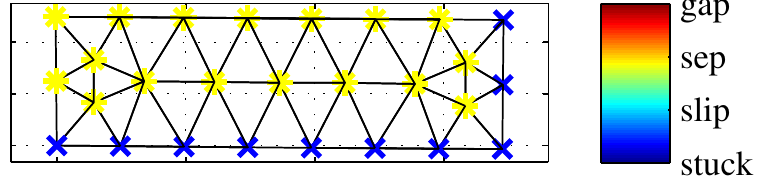}{}{}{.45}{.45}{Contact status for different
modal amplitudes (~(a) point \kreis{3}, (b) point \kreis{4}~)}
\\
The modal properties for the bladed disk are depicted in
\frefs{fig14a}-\frefo{fig14b}.
The qualitative dependency on the energy is similar to the clamped
beam with friction nonlinearity in
\frefs{fig10a}-\frefo{fig10b}.
Apparently, several harmonics have to be considered in the
multi-harmonic analysis in order to achieve asymptotic behavior of
the modal properties due to the highly nonlinear contact
constraints. In this case, more harmonics are required for the
accurate prediction of the modal damping than for the prediction of
the eigenfrequency. Compared to the results in
\fref{fig10b}, the modal damping is less
smooth with respect to the kinetic energy. The reason for this is
that the contact situations change at different energy levels for
each contact node, which can be well ascertained from
\frefs{fig15a}-\frefo{fig16b}. For
small modal amplitudes, only a small portion of the contact area
undergoes stick-slip and partial lift-off. This portion increases
with modal amplitude until the entire set of contact nodes undergo
stick-slip and lift-off phases during the vibration cycle, except of
course for the nodes that were artificially fixed.
\myf[t!]{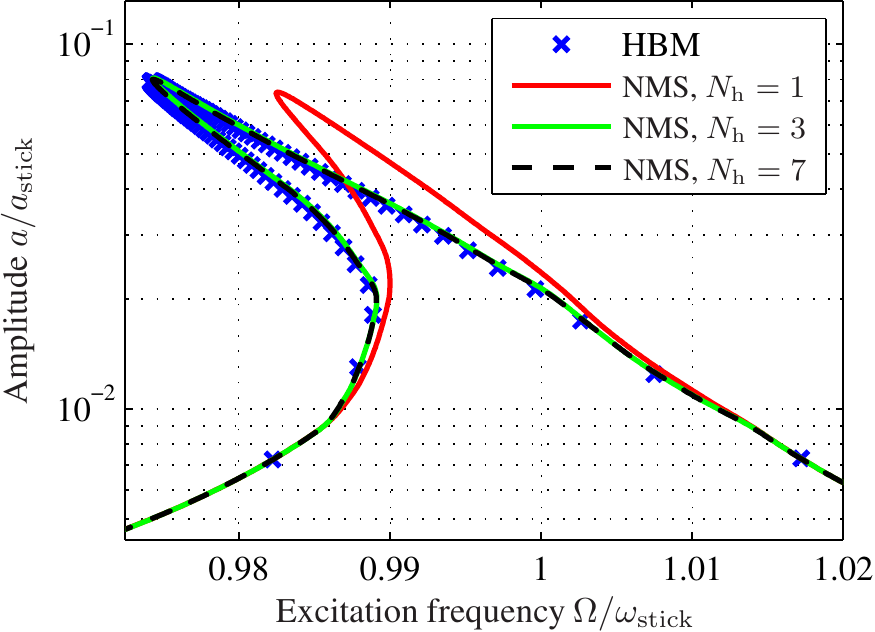}{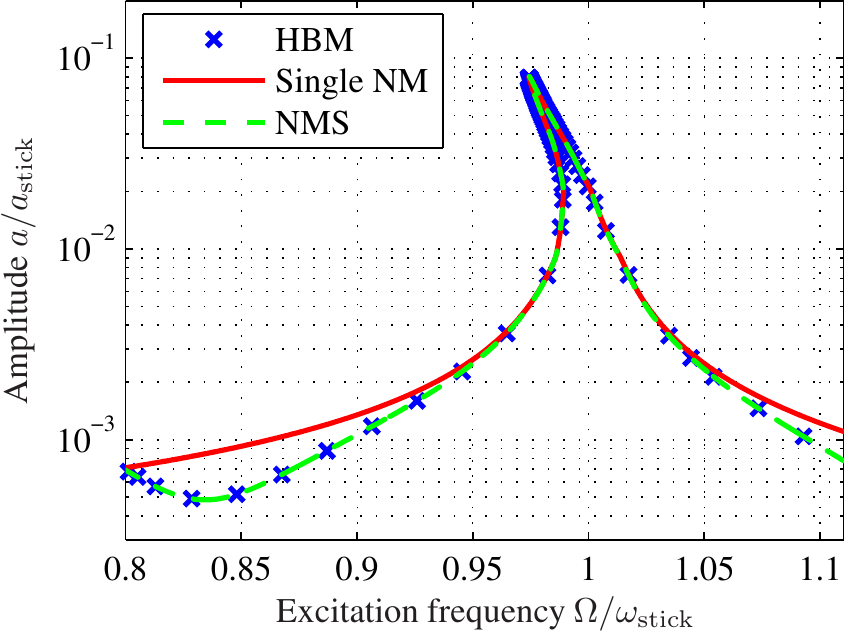}{}{}{.45}{.45}{Accuracy of the forced response
synthesis (~(a) influence of the
number of harmonics \g{\nh} in the modal analysis, (b) influence of the
linearized modes~)}
\\
The accuracy of the forced response synthesis also depends on the
number of harmonics \g{\nh} in the modal analysis, see
\fref{fig17a}. A number of \g{\nh=3} or \g{\nh=7} should
be sufficient to achieve good agreement with the results obtained by
the (multi-term) HBM in this case. From the results in
\fref{fig17b}, the effect of the superposition of the
linearized mode shapes can be deduced. Particularly in the regime
further away from the resonance, the contribution of the linearized
modes becomes more significant and should be accounted for. Note
that the superposition of the linearized modes is a cheap
post-processing calculation and does not significantly increase the
computational expense of the numerical investigations.
\fss[t!]{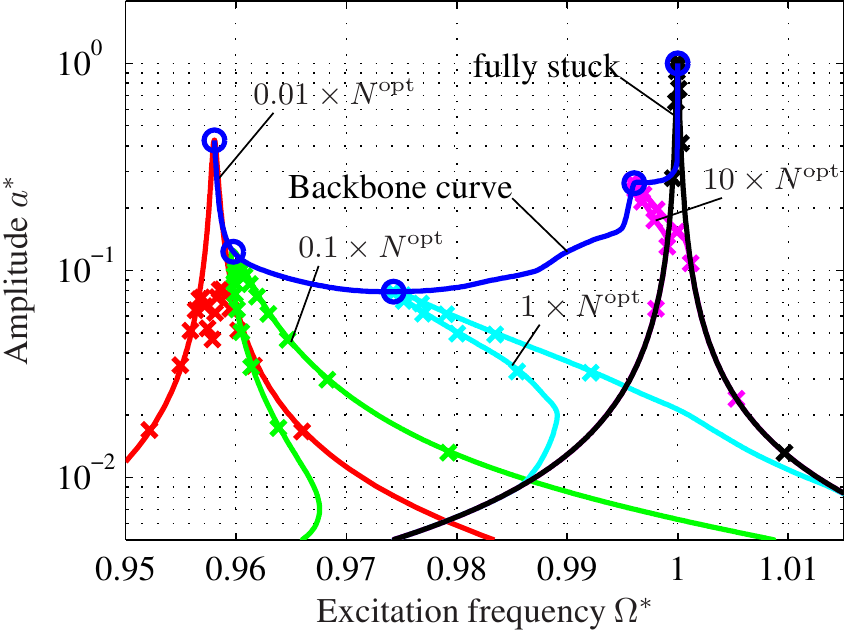}{Forced response of a bladed disk
with shroud contact for varying interlock load}{1.0}
\\
Finally, the forced response was calculated for varying interlock
load \g{N} in the shroud joint, see \fref{fig18}.
Solid lines represent the results obtained by nonlinear modal
synthesis, crosses illustrate the results obtained by the multi-term
HBM with \g{\nh=7}. For large values of \g{N}, the shroud is fully
stuck. For decreasing normal preload, the resonance amplitude is
significantly reduced by means of friction damping up to an optimal
preload value \g{\fnpopt}. The resonance amplitude then increases
again. Below a certain value of the normal preload, the system
exhibits modal interactions. As a consequence, more than one maximum
exists in the forced response. As the assumption of the absence of
internal resonances is no longer valid in this case, the prediction
by the nonlinear modal synthesis (NMS) fails, \cf the results for
\g{N=0.01\fnpopt}. The dynamics of the underlying system cannot be
approximated with only a single nonlinear coordinate anymore. Up to
this regime, however, the prediction of the frequency response
function as well as the backbone curve is in
excellent agreement with the results obtained from the HBM.
\myf[t!]{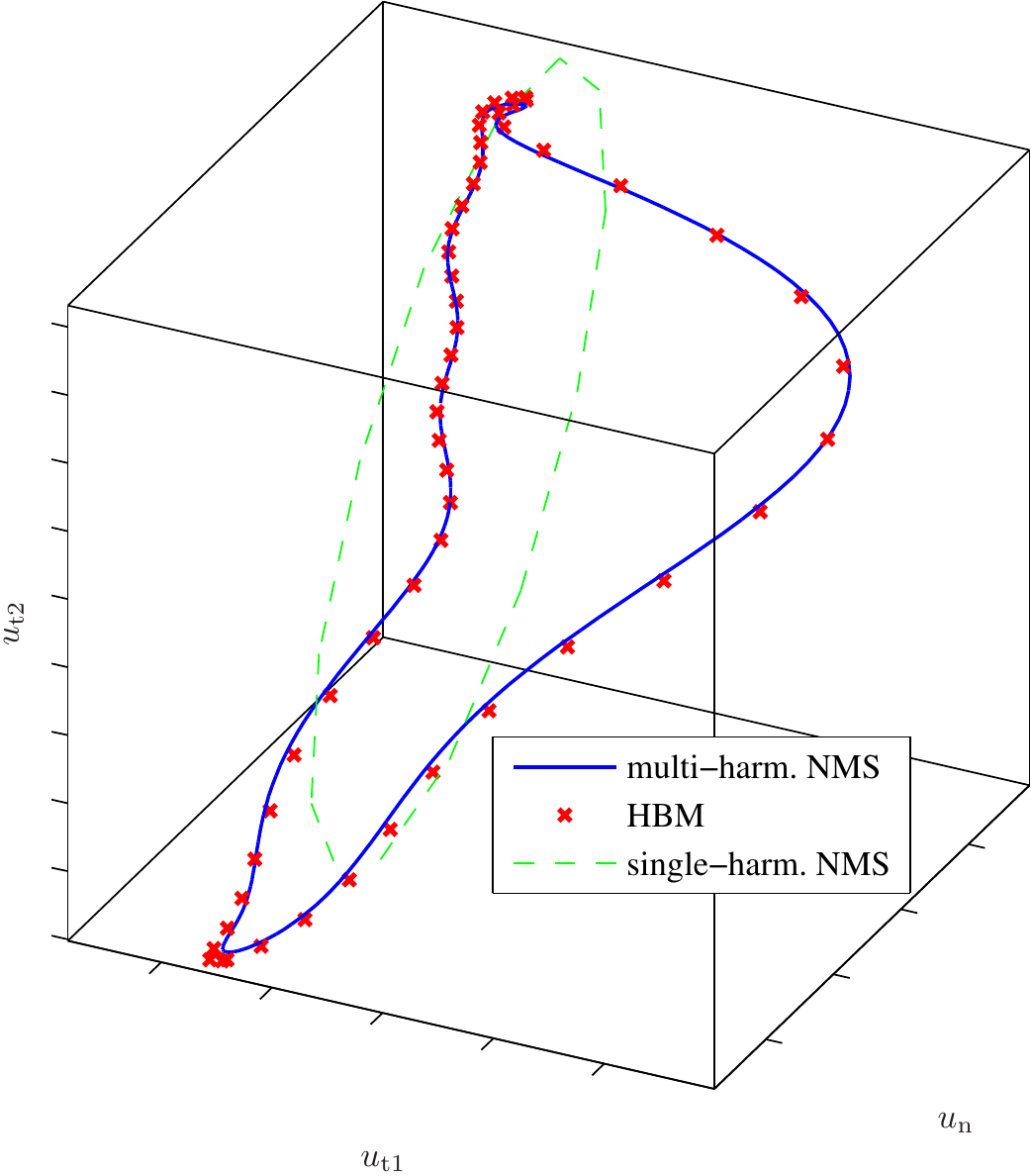}{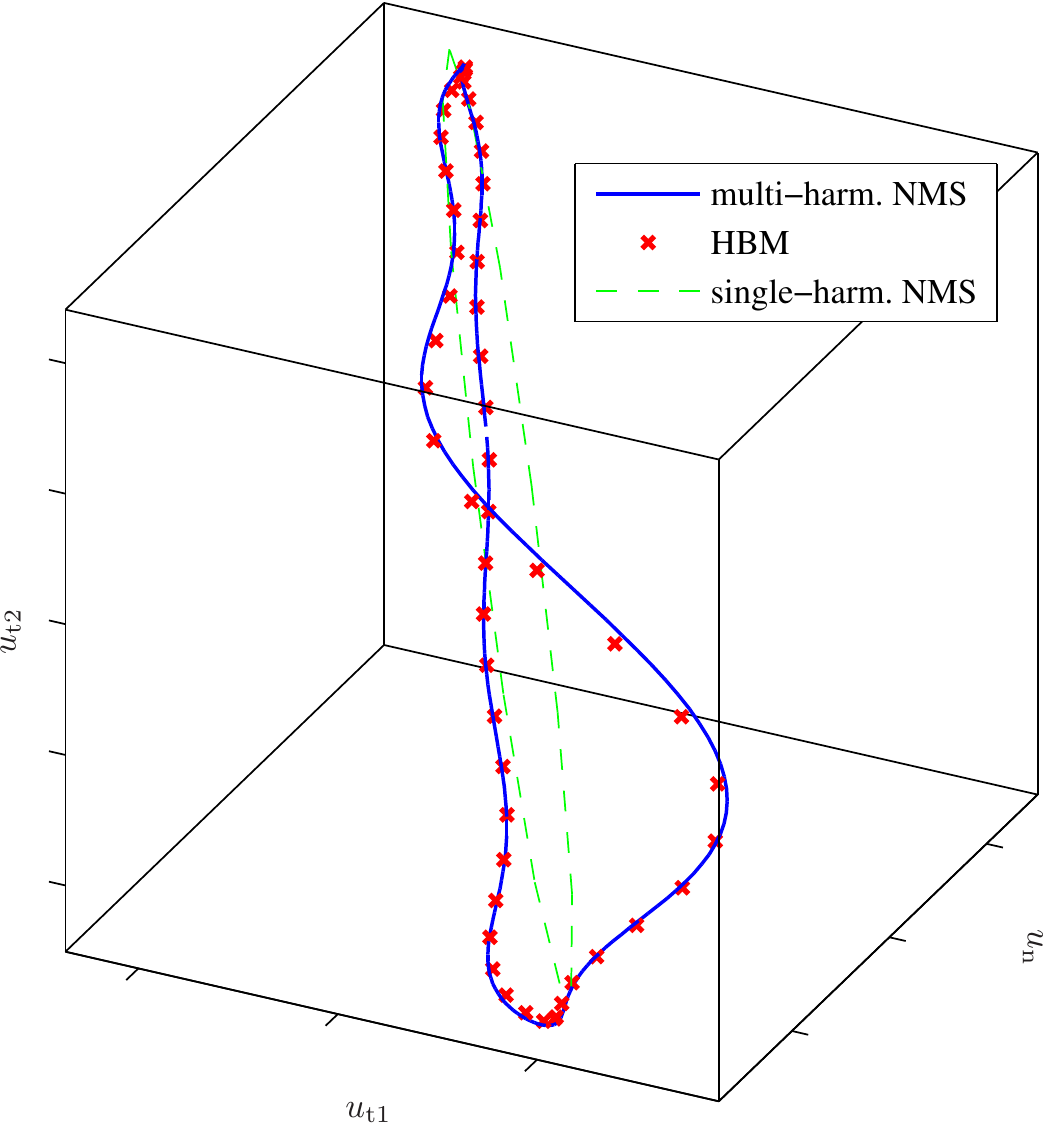}{}{}{.45}{.45}{Orbits of upper
left contact node at resonance (~(a) \g{0.1\times N^{\mathrm{opt}}},
(b) \g{10\times N^{\mathrm{opt}}}~)}
\\
Instead of re-computing the nonlinear modal basis for each value of
the interlock load, the forced response was calculated using the
similarity hypothesis in \eref{preload_excitation_level}. Hence, the
modal properties only had to be calculated once to obtain the
results in \fref{fig18}. Regarding the agreement of
the results, it can be concluded that the similarity hypothesis was
not violated in the considered case.\\
In \frefs{fig19a}-\frefo{fig19b}, the orbits of the
upper left contact node (see \fref{fig13a}) at resonance
are illustrated. It can be ascertained from the results that the
proposed multi-harmonic synthesis significantly increases the
accuracy compared to only considering the fundamental harmonic of
the nonlinear mode.
\tab[b!]{lcc}{Analysis & $N_{\mathrm{dim}}$ & normalized CPU time\\
\hline\hline
FRF (HBM) & $765$ & $1.0$\\
Backbone (HBM) & $766$ & $0.9$\\ \hline
FRF (NMS) & $1$ & $<0.0002$\\
Backbone (NMS) & $1$ & $<0.0002$\\
NMA & $767$ & $2.0$ }{Computational effort for conventional and
proposed methodology}{comp_effort}
\\
In \tref{comp_effort} the computational effort for the forced
response analysis of the bladed disk with shroud contact is listed.
The effort for the nonlinear modal analysis (NMA) of a single mode
is in the order of magnitude of a single frequency response function
(FRF) or backbone curve calculation using the conventional
multi-harmonic balance method (HBM). It should be noted that this
computational effort is strongly related to the number of nonlinear
displacement unknowns
\g{N_{\mathrm{dim}}=(2\nh+1)\nnl=(2\cdot7+1)\cdot 3\cdot 17=765} for
the \g{17} three-dimensional contact elements used in this example,
see \fref{fig13a}. Once the modal properties are known, the
evaluation of the nonlinear modal synthesis (NMS) has almost
negligible computational cost owing to the fact that the number of
unknowns is unity, \ie \g{N_{\mathrm{dim}}=1}.

\section{Conclusions\label{sec:conclusions}}
The recently developed complex nonlinear modal analysis technique
has been refined in this study. With the extensions, it is now
possible to exploit sparsity of the governing algebraic system of
equations, making it particularly attractive for systems featuring
localized nonlinearities. It was also demonstrated that the use of
numerical continuation can facilitate the investigation of modal
interactions with this method. Moreover, it was indicated that the
approach is closely related to the conventional harmonic balance
approaches, so that existing implementations can easily be augmented
with only slight modifications.\\
The resulting nonlinear modal basis was then incorporated into a
novel, very compact ROM based on the single-nonlinear-resonant-mode
theory. Scalar nonlinear equations have been derived for the
calculation of frequency response functions, backbone curves of the
forced response and limit-cycle-oscillations. It was shown that
system parameters, in particular parameters defining the linear
damping and excitation terms can be varied without the need for the
comparatively expensive re-computation of the modal basis. The
proposed technique can thus be employed to facilitate exhaustive
parametric studies on the steady-state vibrations of nonlinear systems.\\
Numerical examples have shown a broad applicability of the overall
methodology. Case studies included large-scale finite element
models, strong and non-smooth, conservative and non-conservative
nonlinearities. In the absence of modal interactions, the synthesis
method showed very good agreement of the multi-modal, multi-harmonic
response with results obtained by conventional methods.\\
Future work on this subject could include the extension of the
nonlinear modal synthesis to transient problems arising \eg in case
of dissipative autonomous systems or in presence of transient
forcing. Bifurcation and stability analyses are considered as
important future developments for the modal analysis technique, in
particular to further investigate nonlinear modal interactions.
Moreover, it would be desirable to extend the ROM to the treatment
of internal resonances which has already been achieved for the
invariant manifold approach. It is conceivable that this could be
accomplished similar to \zo{pierre2006} by an increase of the number
of nonlinear modal amplitudes in modal analysis and synthesis.

\section{Acknowledgements}
The support of Siemens Energy and MTU Aero Engines including the
permission to publish this work is kindly acknowledged. The work
presented in this paper was funded by AG Turbo 2020; Teilvorhaben
4.1.3, FK 0327719A. The responsibility for the content of the
publication rests with the authors.

\begin{appendix}
\section{Dynamic compliance of a system with general structure\label{asec:hgen}}
The dynamic compliance matrix is derived for a system with
invertible mass matrix \g{\mm M} but otherwise general structure.
The dynamic compliance matrix can be computed blockwise for each
harmonic \g{n}. The corresponding dynamic stiffness matrix for the
$n$th harmonic reads,
\e{\mm S_n(\lambda) = \mm K + n\lambda\mm C+ (n\lambda)^2\mm
M\fp}{sngen}
In order to efficiently compute the dynamic compliance matrix \g{\mm
H_n=\mm S_n^{-1}}, a spectral decomposition of the state-space
matrix \g{\mm A} of the system,
\e{\mm A = \matrix{cc}{\mm 0 & -\mm I\\ \mm M^{-1}\mm K & \mm M
^{-1}\mm C}\fk}{matrixa}
is carried out that is defined in analogy to \eref{evlin},
\e{\mm x^{(l)}_k\mm A\mm x^{(r)}_k = \nu_k\,,\,\, \mm x^{(l)}_k\mm
x^{(r)}_k = 1\,,\quad k=1,\cdots,2\ndim\fp}{evagen}
Herein, \g{\mm x^{(l)}_k}, \g{\mm x^{(r)}_k} are left and right
eigenvectors associated with the eigenvalue \g{\nu_k}. Uniqueness of
the eigenvalues is assumed so that non-identical eigenvectors are
orthogonal to \g{\mm A} and to each other. The eigenvectors can be
divided into two blocks of equal dimensions,
\e{\mm x^{(l)}_k = \matrix{cc}{\mm v^{(l)}_k & \mm
w^{(l)}_k}\,,\quad \mm x^{(r)}_k = \vector{\mm v^{(r)}_k \\ \mm
w^{(r)}_k}\fp}{evdiv}
With these definitions and some algebraical manipulations, the
dynamic compliance matrix can finally be identified as
\e{\mm H_n(\lambda) = \suml{k=1}{2\ndim}{\frac{\mm v_k^{(r)}\mm
w_k^{(l)}\mm M^{-1}} {\nu_k+n\lambda}}\,,\quad n =
0,\cdots,\nh\fp}{hinvgen}
Of course, the products \g{\mm w_k^{(l)}\mm M^{-1}} can be carried
out once and for all prior to the nonlinear dynamic analysis.
\end{appendix}





\end{document}